\documentclass[11pt]{article}
\usepackage[a4paper, margin=1in]{geometry}
\usepackage[table]{xcolor}
\usepackage{tcolorbox}
\usepackage{graphicx}
\usepackage{subcaption}
\usepackage[colorlinks,citecolor=blue,urlcolor=blue,pagebackref=true]{hyperref}
\usepackage{amsmath,amsfonts,amssymb,amsthm,euscript,pifont}
\usepackage{mathrsfs}
\usepackage{upgreek}
\usepackage {epigraph}
\usepackage{authblk}
\usepackage{appendix}
\usepackage{multirow}
\usepackage{xspace}
\usepackage{soul}
\usepackage{enumitem}

\usepackage{ulem}

\newtheorem{thm}{Theorem}[section]
\newtheorem{prop}[thm]{Proposition}
\newtheorem{lem}[thm]{Lemma}

\newtheorem{defi}[thm]{Definition}

\newtheorem{remark}[thm]{Remark}

\numberwithin{equation}{section}

\newcommand{\R}{\mathbb{R}}  

\newcommand{\N}{\mathbb{N}}
\newcommand{\E}{\mathbb{E}} 
\newcommand{\hS}{\widehat{\mathscr{S}}}  
\newcommand{\tS}{\widetilde{\mathscr{S}}}  
\newcommand{\hQ}{\widehat{{Q}}}  
\newcommand{\hB}{\widehat{{B}}}

\newcommand{\D}{\Delta} 
\renewcommand{\P}{\mathbb{P}} 
\newcommand{\Pconv}[1]{\xrightarrow[#1]{\P}}
\newcommand{\ASconv}[1]{\xrightarrow[#1]{\text{a.s.}}}
\newcommand{\Dconv}[1]{\xrightarrow[#1]{\mathcal{D}}}

\usepackage{dsfont}
\newcommand{\ind}[1]{\mathds{1}_{#1}}
\newcommand{\zi}{\underline{z}^i}

\DeclareMathOperator{\Cor}{Corr} 
\DeclareMathOperator{\real}{real} 

\DeclareMathOperator{\sgn}{sgn}
\DeclareMathOperator{\DA}{DA} %
\DeclareMathOperator{\law}{law} %

\begin{document}


\title{A hypothesis test for the domain of attraction of a random 
variable}

\author{H\'ector Olivero\thanks{CIMFAV - Instituto de Ingeniería Matemática, Universidad de Valpara\'iso. Chile. Partially supported by FONDECYT POSTDOCTORADO 3180777.}   \and Denis Talay\thanks{INRIA Saclay, France.} }
\date{\today}
 \maketitle

\begin{abstract}

In this work we address the problem of detecting whether a sampled 
probability distribution of a random variable $V$  has infinite first moment. 
This issue
is notably important when the sample results from complex numerical 
simulation methods. For example, such a situation occurs when one 
simulates 
stochastic particle systems with complex and singular McKean-Vlasov
interaction kernels. As stated, the detection problem is ill-posed.
We thus propose and analyze an asymptotic hypothesis test for
independent copies of a given random variable which is supposed to 
belong to an unknown domain of attraction of a stable law. 
The null hypothesis $\mathbf{H_0}$ is: `$X=\sqrt{V}$ is in the domain of 
attraction of the Normal law' and the alternative hypothesis is 
$\mathbf{H_1}$: `$X$ is in the
domain of attraction of a stable law with index smaller than 2'.
Our key observation is that~$X$ cannot have a finite 
second moment when $\mathbf{H_0}$ is rejected (and therefore 
$\mathbf{H_1}$ is accepted).

Surprisingly, we find it useful to derive our test from the statistics
of random processes. More precisely, our hypothesis test is based
on a statistic which is inspired by methodologies to
determine whether a semimartingale has jumps from the observation of
one single path at discrete times.

We justify our test by proving asymptotic properties of discrete time
functionals of Brownian bridges. 

\end{abstract}

\tableofcontents

\section{Introduction}

In this paper we consider situations where one observes a sample of an 
unknown probability distributions and wonder if that probability 
distribution has finite moments.

Let us give an important example of such a situation, which actually 
has motivated our research. A now standard stochastic numerical method
to solve non-linear McKean-Vlasov-Fokker-Planck equations consists in 
simulating stochastic particle systems with McKean-Vlasov interactions
of the type
\begin{equation} \label{particle-system}
X^{i,N}_t =  X^{i,N}_0 
+ \int_0^t \frac{1}{N} \sum_{j=1}^N\beta(X^{i,N}_s,X^{j,N}_s)~ds 
+ \int_0^t \frac{1}{N} \sum_{j=1}^N\gamma(X^{i,N}_s,X^{j,N}_s)~dW^i_s, 
~~1\leq i\leq N,
\end{equation}
where the functions $\beta$ and $\gamma$ from~$\R^d$ to~$\R^d$ are the 
interaction kernels,
the $X^{i,N}_0$'s are independent copies of a random variable~$X_0$
and the $W^i$ are independent Brownian motions. Under appropriate
conditions, the preceding system propagates chaos: The measure-valued
process $\frac{1}{N}\sum_{i=1}^N \delta_{X^{i,N}_\bullet}$ weakly
converges to a unique probability distribution on the space of
continuous functions $\mathcal{C}$ from~$\R^d$ to~$\R^d$ which is the 
probability law of the solution
to the following McKean-Vlasov stochastic differential equation:
\begin{equation} \label{eds:McK-V} 
\begin{cases}
X_t = X_0 
+ \int_0^t \int_{\mathcal{C}} \beta(X_s,\omega_s)~\mu(d\omega)~ds
+ \int_0^t \int_{\mathcal{C}} \gamma(X_s,\omega_s)~\mu(d\omega)~dW_s, \\
\mu~:=~\text{Probability distribution of the process}~(X_t).
\end{cases}
\end{equation}
The time marginal distributions of this
mean-field limit probability measure solve the 
McKean-Vlasov-Fokker-Planck equation under consideration in the sense
of the distributions. See the seminal survey 
by~Sznitman~\cite{Sznitman}. {We emphasize that, within the 
context of
singular interactions, most often it seems inaccessible to prove or 
even intuit accurate tail estimates on the probability 
distributions~$\mu_s$.}

When the interactions kernels $\beta$ and $\gamma$ are so singular than 
no known theorem guarantees the well-posedness of the 
equation~\eqref{eds:McK-V}, numerical simulations may help to study it.
In particular, simulations may give intuition on the (possibly local in 
time) finiteness of the expectations $\int_{\mathcal{C}} 
|\beta(x,\omega_s)|~\mu(d\omega)$ and $\int_{\mathcal{C}} 
|\gamma(x,\omega_s)|~\mu(d\omega)$ with $x\in\R^d$. Of course, this 
property is crucial to get existence of a solution. Similarly, we are 
interested in knowing whether $|\beta(x,\omega_s)|^2$ and 
$|\beta(x,\omega_s)|^2$ are $\mu$-integrable and therefore whether
we can apply a central-limit theorem to get the weak convergence rate
of~$(X^{1,N}_t)$ to $(X_t)$ in terms of the number~$N$ of particles.
See Talay and Tomasevic~\cite{Talay-Tomasevic} for an example of 
numerical simulations within a highly singular context.

Supposing that the system~\eqref{particle-system} is well-posed and 
propagates chaos for small times, considering that the empirical 
measure $\frac{1}{N}\sum_{i=1}^N \delta_{X^{i,N}_s}$ is a good 
approximation of the probability distribution of~$X_s$, the 
numerical strategy the naturally consists in using it to test the 
finiteness of the expectations {or variances} of $|\beta(x,X_s)|$ 
and $|\gamma(x,X_s)|$. 

In~\cite{HAWKINS1997273} Hawkins examines a related question in the 
simpler case of i.i.d. sequences: `Does there exist a test, 
which makes the right decision
with arbitrarily high probability if given sufficient data,
of the hypothesis that a given random variable~$V$ has finite 
expectation?'. 
To address that question, the author introduces the set
$\mathcal{G}$ (resp. $\mathcal{H}$) of densities 
with finite (resp. infinite) means. He also introduces the class~
$\mathcal{T}$ of sequential tests which terminate
in finite time, whatever is the density of the data.
The sets $\mathcal{G}$ and $\mathcal{H}$ would be distinguishable in 
$\mathcal{T}$ if: $\forall\epsilon$ it would exist a test $\phi$
in~$\mathcal{T}$ s.t. 
\begin{equation*}
\begin{cases}
\mathbb{E}_F(\phi)<\epsilon &\text{if}~~ F\in\mathcal{G}, \\
\mathbb{E}_F(\phi)>1-\epsilon &\text{if}~~ F\in\mathcal{H}.
\end{cases}
\end{equation*}
It is proven in~\cite{HAWKINS1997273} that $\mathcal{G}$ and 
$\mathcal{H}$ are 
not distinguishable.

Hawkins' nice result led us to develop a test 
with a more restrictive objective than testing the finiteness of 
$\E|V|$. We actually consider the related question: `How to determine 
from samples of a random variable~$V$ whether its probability 
distribution has heavy tails?'

In our Supplementary Material \cite{SupplementaryMaterial}
we present some numerical experiments which illustrate 
Hawkins' Theorem \cite{HAWKINS1997273} and  
the irrelevance of naive statistical procedures to 
determine the heaviness of the tails of a probability
distribution.

 In this article we construct and analize an asymptotic statistical test under the additional 
assumption that $X:=\sqrt{|V|}$ belongs to some domain of attraction
$\DA(\alpha)$ of a stable law of index $0<\alpha\leq2$
(see Section~\ref{sec:reminders-stable-laws} below for reminders on 
domains of attractions and stable laws).  More precisely, 
 we develop
an hypothesis test for which the null and alternative hypotheses 
respectively are:
\begin{gather*}
\mathbf{H_0}: X\in\DA(2) \\ 
\text{and} \\ 
\mathbf{H_1}: \exists0<\alpha<2,~~X\in\DA(\alpha).
\end{gather*}
\textbf{Our key observation is that $X$ cannot have a finite 
second moment when $\mathbf{H_0}$ is rejected (and therefore 
$\mathbf{H_1}$ is accepted).}

The construction of an effective test is not obvious because, first, the hypothesis test concerns tail properties of the unknown distribution  and, second, when $\alpha$ is close to $2$ the probability density of the symmetric $\alpha$-stable law can hardly be distinguished from the standard Gaussian density on large finite intervals (see for example ~\cite[Figs.2 and 3]{holt1973tables}). Our construction of an effective hypothesis test is original. 
Unexpectedly, it is based on fine 
properties of bivariations of semimartingales and a test for jumps 
which allows one to discriminate between discontinuous stable processes 
and Brownian motions.

Here is our main result.

\begin{thm}
Assume that $X$ belongs to some domain of attraction. Consider and 
i.i.d. sample $X_1,\ldots,X_m$ of $X$, and the statistic
\begin{equation*}
\hS_n^m  = \frac{\sum_{i=1}^{n-1}{\left| \sum_{j=\lfloor 
\frac{m(i-1)}{n} 
\rfloor +1}^{\lfloor \frac{mi}{n} \rfloor } (X_j-\overline{X}_m) 
\right|\left| \sum_{j=\lfloor \frac{mi}{n} \rfloor +1}^{\lfloor 
\frac{m(i+1)}{n} \rfloor } (X_j-\overline{X}_m) 
\right|}}{\sum_{i=1}^{n}{\left|\sum_{j=\lfloor \frac{m(i-1)}{n} 
\rfloor +1}^{\lfloor \frac{mi}{n} \rfloor } (X_j-\overline{X}_m) 
\right|^2} 
},
\end{equation*}
where $\overline{X}_m$ stands for the sample mean. 
Let $z_q$ denote the $q$-quantile of a standard normal random
variable and let $\sigma_\pi^2 := 1 +\frac{4}{\pi}- \frac{20}{\pi^2}$. 
The rejection region 
$$C_{n,m} := \left\{\left| \hS_n^m - \frac{2}{\pi} \right|>z_{1-q/2}
\sqrt{\frac{\sigma_\pi^2}{n}}\right\} $$ 
satisfies: 
\begin{enumerate}
\item $\limsup_{n\to\infty}\limsup_{m\to\infty}\P\left(C_{n,m} |
\mathbf{H_0}\right)\leq q$. 
\item
$\lim_{n\to\infty}\lim_{m\to\infty}\P\left(C_{n,m} | \mathbf{H_1} 
\right)=1$.
\end{enumerate}
\end{thm}

\paragraph*{Plan of the paper:} The plan of the paper is as follows. 

In Section~\ref{sec:reminders-stable-laws} we recall important results
on stable laws and domains of attraction.

In Section~\ref{sec:A-consistent-estimator} we introduce 
the statistic~$\hS_n^m$ on which is based our test. 

In Section~\ref{sec:our-hypothesis-test} we present our hypothesis test 
for domains of attraction of stable laws.

In Section~\ref{sec:consistency} we prove the consistency of the 
statistic~$\hS_n^m$. 

In Section~\ref{sec:CLT} we prove a Central Limit theorem for
$\hS_n^m$.

In Section~\ref{sec:num-exp}  we discuss some numerical experiments.

In our conclusion we comment on the performance and limitations of our 
test.

Finally, in the Appendix \ref{sec:appendix} we prove few intermediate 
technical results. 

\section{A few reminders on stable laws and domains of attraction}
\label{sec:reminders-stable-laws}

In this section we gather the few results on stable laws and domains of 
attraction which we need in the sequel. For their proofs and further 
information, see e.g. 
Feller~\cite[Sec.8,Chap.9]{feller1971introduction}, Embrechts 
et al.\cite[Sec.2, Chap.2] {embrechts2013modelling} 
or Whitt~\cite[Sec.5, Chap.4]{whitt2002stochastic}. We also reformulate 
a standard functional limit theorem under a form which prepares our 
hypothesis test to reject or accept~$\mathbf{H_0}$.

Let $X,X_1,X_2,\ldots$ be a sequence of non-degenerate i.i.d. random
variables and let 
$$ S_m := \sum_{j=1}^{m}{X_j},\;\; m=1,2,\ldots $$ 
One says that $X$, or the law of $X$, belongs to the 
domain of attraction of a given law~$\mathcal{L}$
if there exist centering constants $\mu_m$ and positive 
normalizing constants $c_m$ such that $(S_m-\mu_m)/c_m$ converges in 
distribution to~$\mathcal{L}$. 

The only probability laws which have a
non-empty domain of attraction are the stable laws.

The probability distributions of stable laws are fully characterized by 
the following theorem (see e.g. 
Feller~\cite[Thm.1a,Sec.8,Chap.9]{feller1971introduction}). 
We recall that a function $H:\R_+\to\R$ is said to be slowly varying at 
infinity if $\lim_{s\to\infty} H(sx)/H(s)=1$ for any $x>0$. The 
logarithm is an example of such a function. 

\begin{thm} \label{theo:pdf-stable-law}
Let $F$ be a probability distribution function.
For $R>0$ denote the truncated second moment of $F$ by
$$ U(R) := \int_{-R}^{R}x^2F(dx). $$

\begin{enumerate}[label=(\roman*)]
\item The probability distribution $F$ belongs to the domain of 
attraction $\DA(2)$ of the Gaussian distribution if and only 
if the function~$U$ is slowly varying at infinity.
\item It belongs to some other domain of attraction~$\DA(\alpha)$ with 
$0<\alpha<2$ if there exist a slowly varying at infinity function
$H:\R_+\to\R$ and positive numbers $p$ and $q$ such that
\begin{equation*}
\begin{cases}
1-F(x)+F(-x)\sim \frac{2-\alpha}{\alpha}~x^{-\alpha}~H(x),~~~ 
x\to\infty, \\
\lim_{x\to\infty}\frac{1-F(x)}{1-F(x)+F(-x)}= p, 
\;\;\;\lim_{x\to\infty}\frac{F(-x)}{1-F(x)+F(-x)}= q.
\end{cases}
\end{equation*}
\end{enumerate}
\end{thm}

The parameter $\alpha$ is called the characteristic exponent or the 
index of the stable law. Given~$0<\alpha\leq 2$, every stable law
with index~$\alpha$ is called an $\alpha$-stable distribution
and denoted by $\mathcal{L}_\alpha$. 

For specific values of~$\alpha$ the stable law~$\mathcal{L}_\alpha$
is simple. For $\alpha=2$,  
$\mathcal{L}_\alpha$ is Gaussian (and possibly degenerate). For~$\alpha=1$,
$\mathcal{L}_\alpha$ is Cauchy. For $\alpha=1/2$, $\mathcal{L}_\alpha$ is L\'evy. 
For $\alpha<2$, $\mathcal{L}_\alpha$ has a density.

The preceding theorem implies the following moment properties
(see e.g. Embrechts et al.~\cite[Cor.2.2.10]{embrechts2013modelling}).

\begin{thm} \label{theo:moments-stable-law}
If the random variable~$X$ belongs to~$D(\alpha)$ then
\begin{align*}
\E(|X|^\delta) < \infty~~~\text{for}~~~0<\delta<\alpha, \\
\E(|X|^\delta) = \infty~~~\text{for}~~\delta>\alpha~~~\text{and}~~~
\alpha<2.
\end{align*}
In particular, $\E(X^2)=\infty$ for $\alpha<2$.
\end{thm}   

\begin{remark} 
In view of Theorem~\ref{theo:pdf-stable-law}, 
if $\E(X^2)<\infty$ then $X$ belongs to~$\DA(2)$. The converse is not 
true.
For example, $X=\sqrt{|Y|}$ with~$Y$ Cauchy belongs to~$\DA(2)$ and has
an infinite second moment. This observation combined with 
Theorem~\ref{theo:moments-stable-law} explains our choice of the null 
hypothesis.
\end{remark}

In the Introduction we explained why  
often one cannot get a priori theoretical tail estimates 
on the probability distribution of particles with complex
dynamics. 
Consequently, even if one a priori knows 
that this probability distribution belongs to some domain of 
attraction, to
determine the parameter~$\alpha$ is a hard statistical question.
That was a motivation to build a specific hypothesis test which
is based on the following limit theorems.

We start with Donsker's invariance principle for i.i.d. random 
variables. See e.g. Whitt~\cite[Thms 4.5.2-4.5.3]{whitt2002stochastic}.

\begin{thm} \label{thm:donsker} 
Let $X,X_1,X_2,\ldots$ be a sequence of non-degenerate i.i.d. random 
variables such that $X\in\DA(\alpha)$. Then there exist centering 
constants $\mu_m$ and normalizing constants $c_m$ such that
$$ L^m := \left(\frac{S_{\lfloor m t \rfloor } - \mu_m t}{c_m},t\geq0 
\right) 
\Dconv{m\to\infty} L,$$
where $L$ is a standard $\alpha$-stable L\'evy process if $\alpha<2$,
whereas for $\alpha=2$~$L$ is a Brownian motion. Here, one considers
weak convergence in the Skorohod space endowed with $J_1$ topology.
\end{thm}

\begin{figure}[tb!] 
\centering
\begin{subfigure}[t]{0.4\linewidth}
\centering
\includegraphics[width=\linewidth]
{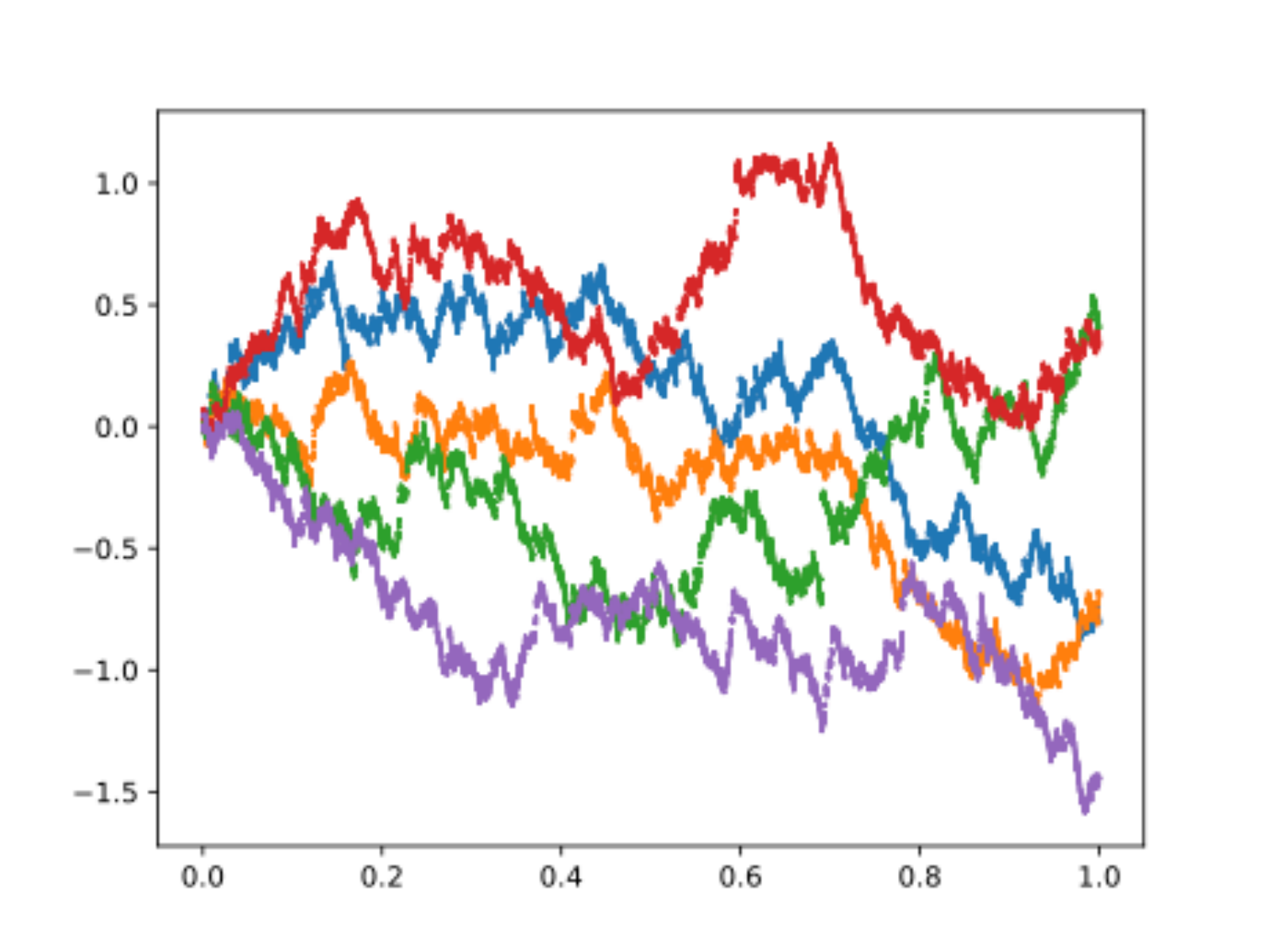}
\caption{\fontsize{9}{11}\selectfont $r=0.2$ 
\label{figure:brownian}}
    \end{subfigure}
    \quad   
\begin{subfigure}[t]{0.4\linewidth}
\centering        
\includegraphics[width=\linewidth]
{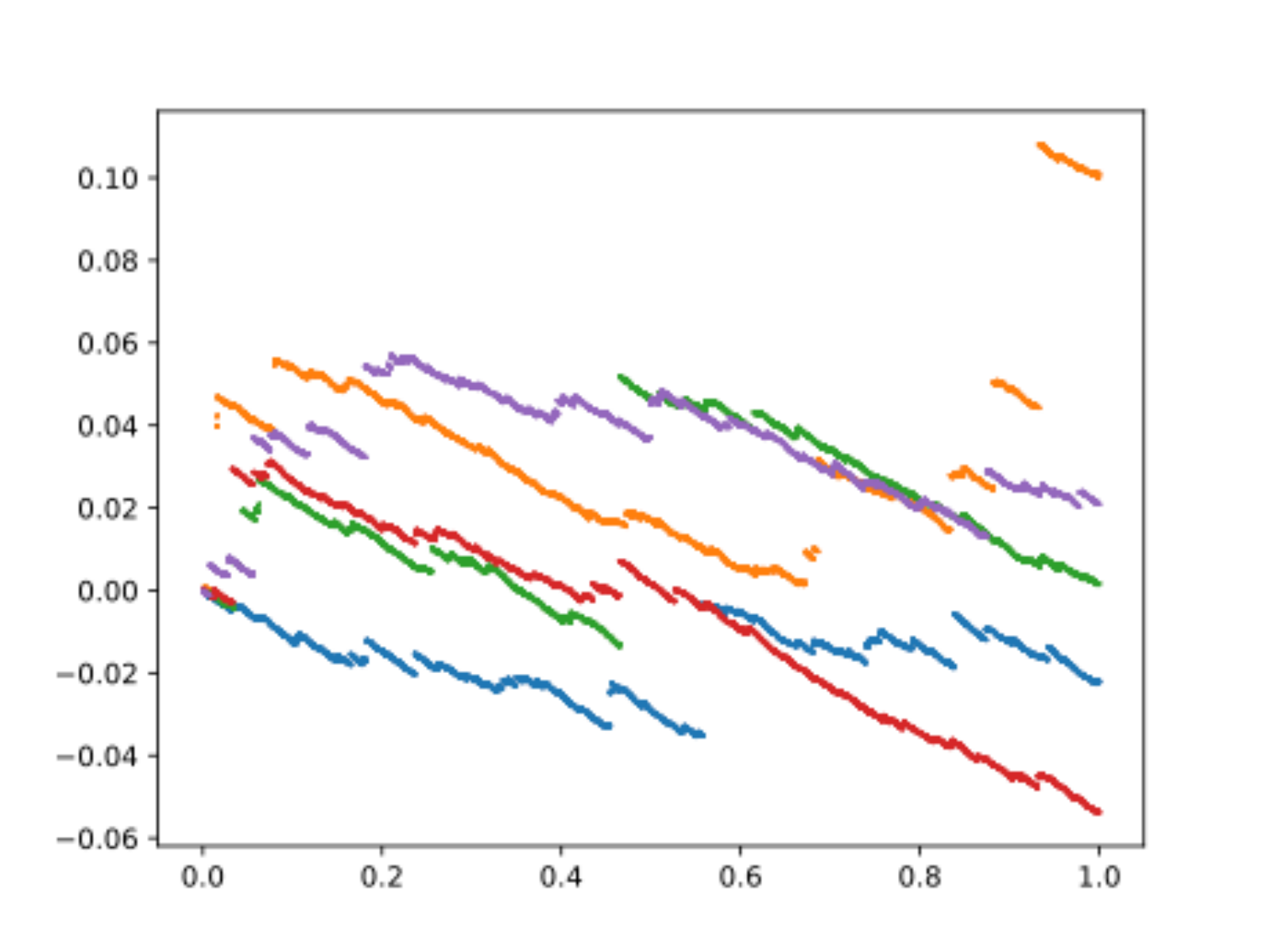}
\caption{\fontsize{9}{11}\selectfont $r=0.8$  
\label{figure:levy}}
\end{subfigure}
\caption{\fontsize{9}{11}\selectfont Trajectories of $L^m$ when the 
subjacent random variable $X\sim |G|^{-r}$, where $G\sim 
\mathcal{N}(0,1)$. In the left column $r=0.2$, therefore the limit 
of $L^m$ is a Brownian motion and the trajectories of $L^m$ 
seem to be continuous. In the right column,
$r=0.8$, hence the limit of $L^m$ is a L\'evy process and the 
trajectories of $L^m$ seem to have jumps. 
\label{figure:donsker-trajectories}}
\end{figure}

For $\alpha<2$, the trajectories of $\alpha$-stable processes are a.s. 
discontinuous, whereas for $\alpha=2$ the trajectories of the
Brownian motion are a.s. continuous. In addition, one expects 
that for $m$ large enough, the trajectories of 
$(S_{\lfloor m t \rfloor } - \mu_m t)/c_m$ 
resemble the trajectories of the limit process (see 
Fig.~\ref{figure:donsker-trajectories}).
Therefore, testing for jumps in the trajectories  
of $(S_{\lfloor m t \rfloor} - \mu_m t)/c_m$ should allow to 
discriminate between $X\in\DA(2)$ and $X\in\DA(\alpha)$. 
{This is illustrated by Fig.~\ref{figure:donsker-trajectories}.
Simulations have been run with $X=|G|^{-r}$, with $G$ a standard normal random variable and  $m = 10000$. For $r=0.2$ we are under 
$\mathbf{H_0}$ and we can compute the explicit value of~$\E(X)$ and 
$\E(X^2)$. For $r=0.8$ we know that $|G|^{-r}$ is in the normal domain 
of attraction of stable distribution with index $1/r$. Hence, the 
normalizing constant is $c_m=m^{1/r}$.} Indeed, respectively denoting 
by $\Phi$ and $\phi$ the 
cumulative distribution function and the density function of a standard 
normal distribution we have
\begin{align*}
	 \P(|G|^{-r}>x) 	&=  2\Phi\left(\frac{1}{x^{1/r}}\right)-1.   
\end{align*}
From L'Hôpital's Theorem it follows that
\begin{align*}
\lim_{x\to\infty} \frac{ \P(|G|^{-r}>x)}{1/x^{1/r}} 	
			=  \lim_{x\to\infty} 
			\frac{2\Phi\left({1}/{x^{1/r}}\right)-1}{{1}/{x^{1/r}}}  
		  	=    \lim_{x\to\infty} 
		  	\frac{2\phi\left({1}/{x^{1/r}}\right)
		  	\left({1}/{x^{1/r}}\right)'}
		  	{\left({1}/{x^{1/r}}\right)'}	 
			=  \sqrt{\frac{2}{\pi}}.
\end{align*}
We thus deduce the desired result 
from~\cite[Thm.4.5.2]{whitt2002stochastic}.

This heuristic approach has a severe drawback: we do 
not know the values of~$\mu_m$ and $c_m$, in particular because 
we do not know~$\alpha$. Nevertheless, as  Proposition  
\ref{cor:convergenceToBridge} below shows that one can use the Mapping 
Theorem to bypass the fact that $\mu_m$ is unknown.

\begin{remark}\label{rem:donsker-under-weak-dependency} The conclusion of 
Theorem \ref{thm:donsker} 
still holds true  in some cases when the $X_i$'s are weakly 
dependent. Sufficient conditions to respectively obtain Brownian or stable L\'evy limits can be found in Shao~\cite{Qiman:1993aa} or Tyran-Kamińska~\cite{TYRANKAMINSKA20101629}.
\end{remark}

The following lemma will allow us to reformulate 
Theorem~\ref{thm:donsker} as a statement where~$\mu_m$ does not need
to be known.

\begin{lem} \label{lem:continuity-of-Psi-Skokokhod-Topology}
Denote by $D[0,1]$ the space of {\it cadlag} functions defined 
on $[0,1]$ and consider 
the mapping $\Psi:D[0,1]\to D[0,1]$ defined as
$$\Psi(x)(t) := x(t)-tx(1).$$
The mapping $\Psi$ is continuous for the the $J_1$-Skorokhod 
topology.
\end{lem}

\begin{proof}
The addition is continuous in the $J_1$-Skorokhod topology if one of the summands is continuous. (See e.g. Jacod and Shiryaev \cite[Prop.1.23,Sec.1b,Chap.VI]{jacod2013limit}).
\end{proof}
 
From the previous lemma and the Mapping Theorem we deduce:
\begin{prop}\label{cor:convergenceToBridge} 
Let 
\begin{equation}\label{eq:defZm}
Z^m:= \left(\frac{S_{\lfloor m t \rfloor } - tS_m}{c_m},t\geq0 
\right).
\end{equation}
If $X\in\DA(\alpha)$, then
$$Z^m \Dconv{m\to\infty} Z:=\Psi(L),$$
where $L$ is a standard $\alpha$-stable L\'evy process for $\alpha<2$
(and therefore $\Psi(L)$ is a discontinuous process), 
whereas for $\alpha=2$~$L$ is a Brownian motion (and therefore
$\Psi(L)$ is a Brownian bridge).
\end{prop}

To conclude, testing for jumps the trajectories of the limit process
of~$Z^m$ defined in~\eqref{eq:defZm} should allow one to discriminate 
between 
$X\in\DA(2)$ and $X\in\DA(\alpha),~0<\alpha<2$.
It now remains to construct an hypothesis test for the continuity 
of~$\Psi(L)$ which does not suppose that~$c_m$ is known.

\section{A bivariation based statistic}
\label{sec:A-consistent-estimator}

Determining whether a stochastic process $Y$ is continuous or not from
the observation of one single path at discrete times, is an important 
modelling issue in many fields, notably in economics and financial 
mathematics. It has been addressed by several authors in the last 
years. See for example  Ait-Sahalia and Jacod~\cite{Ait-Sahalia:2009aa} 
and the references therein. 

Our situation is somehow different since we are not observing a 
trajectory of a given process. We rather are constructing one discrete 
time path by means of a 
normalization procedure of our (observed or simulated) data and this 
constructed trajectory is always 
discontinuous. We therefore aim to construct a statistic whose
asymptotic properties will allow us to apply detection of jumps of 
semimartingale methods.

Following Barndorff-Nielsen and Shephard 
\cite{barndorff2006econometrics}, for any stochastic 
process $(Y_t)_{0\leq t\leq 1}$ we set $\D_i^n Y =  
Y_{i/n}-Y_{(i-1)/n}$ and consider the realized bivariation, the 
realized quadratic variation and the normalized realized bivariation respectively
defined by 
\begin{equation} \label{def:S-hat-n}
\hB(Y,n):=\sum_{i=1}^{n-1} |\D_i^n Y||\D_{i+1}^n Y|,\quad 
\hQ(Y,n):=\sum_{i=1}^{n} |\D_i^n Y|^2, \quad \hS_n(Y):= 
\frac{\hB(Y,n)}{\hQ(Y,n)}.
\end{equation}

In \cite{barndorff2006econometrics}, the authors consider 
processes~$(Y_t)$ of the form
\begin{equation}\label{eq:BSMJ-BNS}
Y_t = Y_0 + \int_{0}^{t}a_s~ds + \int_{0}^{t}\sigma_s~dW_s 
+ \sum_{j=1}^{N_t}c_j,
\end{equation}
where $W$ is a Brownian motion and $(N_t)$ is a counting process.  
They assume that $a$ and $\sigma$ are c\`adl\`ag processes,
$c_j$ are non-zero random variable, $\sigma$ is pathwise bounded away 
from zero and the joint process~$(a,\sigma)$ is independent of the 
Brownian motion~$(W_t)$. Barndorff-Nielsen and Shephard  provided a 
test to decide
`$\mathcal{H}_0$: $(N_t)\equiv 0$' against `$\mathcal{H}_1$: 
$(N_t)\not\equiv 0$'.

Barndorff-Nielsen and Shephard's test is based on the statistic
$$ \frac{1}{\sqrt{n}}
~\left(\frac{\frac{2}{\pi}\hB(Y,n)}{\hQ(Y,n)}-1\right)
~\frac{\int_0^t \sigma_s^2~ds}{\sqrt{\int_0^t \sigma_s^4~ds}} $$
or to a similar statistic depending on $\hB(Y,n)$, $\hQ(Y,n)$
and $\int_0^t \sigma_s^4~ds$ only.

Barndorff-Nielsen and Shephard's  analysis of 
the test (see~\cite[Thm.1]{barndorff2006econometrics}) does not apply 
in our context since
Brownian bridges
are not of the type~\eqref{eq:BSMJ-BNS} with $(\sigma_s)$ satisfying 
the requested constraints.
 
The non-degeneracy of the integrands is also needed by Ait-Sahalia and 
Jacod (see Assumption 1-(f) in~\cite[p.187]{Ait-Sahalia:2009aa}).
In addition, they use the p-variation of semimartingales with $p>2$ to 
test continuous paths against jumps. We cannot follow the same way 
since we test `Brownian bridge' against `a transformed alpha-stable 
process' (recall~Proposition~\ref{cor:convergenceToBridge}) and
the semimartingale representation of a Brownian bridge 
implies that the p-variation is asymptotically infinite in both cases.

We therefore need to introduce a new statistic which is adapted to
our specific situation.

\begin{defi} For any $n\in\N$ consider the real-valued functional 
$\hS_n$ defined as follows: for $z\in D[0,1]$,
\begin{equation}\label{eq:functional-hS}
\hS_n(z) := \frac{\sum_{i=1}^{ 
n-1}{|z(i/n)-z((i-1)/n)||z((i+1)/n)-z(i/n)|}}
{\sum_{i=1}^{n}{|z(i/n)-z((i-1)/n)|^2} }.
\end{equation}
Let $m\in\N$, given an i.i.d. sample of $X:X_1,\ldots,X_m$ and
the corresponding process $Z^m$ as in~\eqref{eq:defZm} we define our 
statistic $\hS_n^m$ as the normalized bivariation of~$Z^m$, that is,
\begin{equation}\label{eq:def-Statistic}
\hS_n^m:=\hS_n(Z^m).
\end{equation}
\end{defi}

The next result provides a key property of the statistic~$\hS_n^m$.

\begin{prop}\label{prop:formula-for-hSn}
The statistic $\hS_n^m$ is scale-free and satisfies
\begin{equation*}
 \hS_n^m  = \frac{\sum_{i=1}^{n-1}{\left| \sum_{j=\lfloor \frac{m(i-1)}{n} 
 \rfloor +1}^{\lfloor \frac{mi}{n} \rfloor } (X_j-\overline{X}_m) 
 \right|\left| \sum_{j=\lfloor \frac{mi}{n} \rfloor +1}^{\lfloor 
 \frac{m(i+1)}{n} \rfloor } (X_j-\overline{X}_m) 
 \right|}}{\sum_{i=1}^{n}{\left|\sum_{j=\lfloor \frac{m(i-1)}{n} 
 \rfloor +1}^{\lfloor \frac{mi}{n} \rfloor } (X_j-\overline{X}_m) 
 \right|^2} }.
\end{equation*}
\end{prop}

\begin{proof} 
Observe that
\begin{align*}
\hS_n^m &:= \hS_n(Z^m)\\
 &=     \frac{\sum_{i=1}^{ 
 n-1}{|Z_{i/n}^m-Z_{(i-1)/n}^m|
 |Z_{(i+1)/n}^m-Z_{i/n}^m|}}{\sum_{i=1}^{n}{|Z_{i/n}^m-Z_{(i-1)/n}^m|^2}
  }\\
  &= \frac{\sum_{i=1}^{ n-1}{\left|\frac{S_{\lfloor \frac{mi}{n} 
  \rfloor } - \frac{i}{n}S_m}{c_m} -\frac{S_{\lfloor \frac{m(i-1)}{n} 
  \rfloor } - \frac{(i-1)}{n}S_m}{c_m}  \right|\left|\frac{S_{\lfloor 
  \frac{m(i+1)}{n} \rfloor } - \frac{(i+1)}{n}S_m}{c_m} 
  -\frac{S_{\lfloor \frac{mi}{n} \rfloor } - \frac{i}{n}S_m}{c_m}  
  \right|}}{\sum_{i=1}^{n}{\left|\frac{S_{\lfloor \frac{mi}{n} \rfloor 
  } - \frac{i}{n}S_m}{c_m}-\frac{S_{\lfloor \frac{m(i-1)}{n} \rfloor } 
  - \frac{i-1}{n}S_m}{c_m}  \right|^2} } \\
  &=\frac{\sum_{i=1}^{n-1}{\left| S_{\lfloor \frac{mi}{n} \rfloor } 
  -S_{\lfloor \frac{m(i-1)}{n} \rfloor } - \frac{1}{n}S_m \right|\left| 
  S_{\lfloor \frac{m(i+1)}{n} \rfloor } -S_{\lfloor \frac{mi}{n} 
  \rfloor } - \frac{1}{n}S_m \right|}}{\sum_{i=1}^{n}{\left|S_{\lfloor 
  \frac{mi}{n} \rfloor } -S_{\lfloor \frac{m(i-1)}{n} \rfloor } - 
  \frac{1}{n}S_m  \right|^2} }\\
 &=\frac{\sum_{i=1}^{n-1}{\left| \sum_{j=\lfloor \frac{m(i-1)}{n} 
 \rfloor +1}^{\lfloor \frac{mi}{n} \rfloor } (X_j-\overline{X}_m) 
 \right|\left| \sum_{j=\lfloor \frac{mi}{n} \rfloor +1}^{\lfloor 
 \frac{m(i+1)}{n} \rfloor } (X_j-\overline{X}_m) 
 \right|}}{\sum_{i=1}^{n}{\left|\sum_{j=\lfloor \frac{m(i-1)}{n} 
 \rfloor +1}^{\lfloor \frac{mi}{n} \rfloor } (X_j-\overline{X}_m) 
 \right|^2} },
\end{align*}
which is the desired result.
\end{proof}

\begin{remark} Notice that computing $\hS_n^m$ only needs the values
of the sample. In particular, one does not need the unknown
centering and normalizing factors $\mu_m$ and $c_m$~in 
Theorem~\ref{thm:donsker}. 
\end{remark}

\section{Our hypothesis test for domains of attraction of stable laws} 
\label{sec:our-hypothesis-test}

In this section we present our hypothesis test for $\mathbf{H_0}$
against $\mathbf{H_1}$. It is based on the following 
consistency property and Central 
Limit Theorem for the statistic~$\hS_n^m$.

The consistency of $\hS_n^m$ is provided by the following proposition
whose lenghthy proof is postponed to Section~\ref{sec:consistency}.

\begin{prop}\label{prop:conv-in-P-estimator} 
For any $\epsilon>0$ one has
\begin{equation}\label{eq:conv-in-prob-statistic-Donsker-approx}
\lim_{n\to\infty}\lim_{m\to\infty}\P\left(|\hS_n^m-\kappa|
\leq\epsilon\right) =1, 
\end{equation} with 
$$\kappa :=\left\{\begin{array}{ccl}  
\frac{2}{\pi}& \text{when} &X\in \DA(2), \\
0& \text{when} & X\in \DA(\alpha),\;\alpha<2. 
\end{array}\right.$$ 
\end{prop}

The next proposition provides a Central Limit Theorem for the 
statistic $\hS_n^m$. We postpone its proof to Section~\ref{sec:CLT}.

\begin{prop} \label{prop:central-limit-theorem-statistic}
If the subjacent random variable $X$ belongs to the domain
of attraction of the normal law,  for any  bounded and continuous function  
$\varphi:\R\to\R$ we have
$$ \lim_{n\to\infty}\lim_{m\to\infty}
\E\left[\varphi\left(\frac{\sqrt{n}}{\sigma_\pi}\left(\hS_n^m
-\frac{2}{\pi}\right)\right)\right]
= \E\left[\varphi(\mathcal{N})\right],$$
where $\sigma_\pi^2 =1 
+\frac{4}{\pi}- \frac{20}{\pi^2} $ {and $\mathcal{N}$ is
a standard Gaussian random variable}.
\end{prop}

Propositions~\ref{prop:conv-in-P-estimator} 
and~\ref{prop:central-limit-theorem-statistic}
allow us to propose and analyse the following hypothesis test:

\begin{thm}\label{thm:TEST}  
Assume that $X$ belongs to some domain of 
attraction. Consider and i.i.d. sample $X_1,\ldots,X_m$ of a r.v.~$X$. 
We consider the test hypotheses

\noindent $\mathbf{H_0}$: $X\in\DA(2)$\\ and \\ $\mathbf{H_1}$: 
$\exists 0<\alpha<2,\;X\in\DA(\alpha)$. \\

Let $z_q$ denote the $q$-quantile of a standard normal random
variable and $\sigma_\pi^2 =1 
+\frac{4}{\pi}- \frac{20}{\pi^2} $. The rejection region 
$$ C_{n,m} := \left\{\left|\hS_n^m-\frac{2}{\pi}\right|>z_{1-q/2}
{\frac{\sigma_\pi}{\sqrt{n}}}\right\} 
{
= \left\{{\frac{\sqrt{n}}{\sigma_\pi}}~(\hS_n^m-\frac{2}{\pi}) < z_{q/2}
\right\} \bigcup  \left\{{\frac{\sqrt{n}}{\sigma_\pi}}~(\hS_n^m-\frac{2}{\pi}) > z_{1-q/2}
\right\} 
} 
$$ 
satisfies: 
\begin{enumerate}
\item $\limsup_{n\to\infty}\limsup_{m\to\infty}\P\left(C_{n,m} |
\mathbf{H_0}\right)\leq q$. 
\item
$\lim_{n\to\infty}\lim_{m\to\infty}\P\left(C_{n,m} | 
\mathbf{H_1}\right)=1$.
\end{enumerate}
\end{thm}

\begin{proof}
To prove the first claim, we fix $\delta>0$ and consider the function 
$\psi_\delta:\R \to\R$ given by:
$$
\psi_\delta(x)=\left\{
\begin{array}{ccl}
1& \text{for}& x<z_{q/2}, \\
1-\frac{1}{\delta}(x-z_{q/2}) &\text{for}& z_{q/2} \leq x \leq z_{q/2}+\delta, \\
0 &\text{for} & z_{q/2}+\delta < x < z_{1-q/2}-\delta, \\
\frac{1}{\delta}(x-z_{1-q/2}) &\text{for}& z_{1-q/2} \leq x \leq z_{1-q/2}, \\
1 &\text{for} & z_{1-q/2} < x. \\
\end{array}
\right.
$$
Therefore,
$$ \P\left(C_{n,m} | \mathbf{H_0} \right) \leq 
\E\left(\psi_\delta(\hS_n^m )|\mathbf{H_0}\right). $$
In view of Proposition~\ref{prop:central-limit-theorem-statistic} we 
obtain 
$$ \limsup_{n\to\infty}\limsup_{m\to\infty}\P\left(C_{n,m} | 
\mathbf{H_0}\right) 
\leq \E\left[\psi_\delta(\mathcal{N})\right],$$
where $\mathcal{N}$ is a standard Gaussian random variable. Thanks to  
the dominated convergence theorem we have
$$\lim_{\delta\to0}  \E\left[\psi_\delta(\mathcal{N})\right] = 
\P\left(\mathcal{N}< z_{q/2}\right)+\P\left(\mathcal{N}> z_{1-q/2}\right)=q. $$
We thus can conclude that
$$ \limsup_{n\to\infty}\limsup_{m\to\infty}\P\left(C_{n,m} | 
 \mathbf{H_0}\right) \leq q. $$
As for the second claim, it suffices to use
Proposition~\ref{prop:conv-in-P-estimator} since under $\mathbf{H_1}$ 
it holds that 
$$ \forall \epsilon>0, 
~~\lim_{n\to\infty}\lim_{m\to\infty}
\P\left({\hS_n^m} \leq \epsilon\right)=1. $$ 
\end{proof}

\begin{remark}
{\bf Other tests in the literature related to ours:} In 
\cite{jurevckova2001class} Jure{\v{c}}ková and Picek develop a 
statistical test for the heaviness of the tail of a distribution 
function $F$ assuming that $F$ is absolutely continuous and strictly 
increasing on the set $\{x:F(x)>0\}$. In their case, for $m_0$ given, 
the null hypothesis is 
$$
H_{m_0}: \;x^{m_0}(1-F(x))\geq1,\;\;\forall x>x_0 \text{ for some } 
x_0\geq0,
$$
whereas the alternative hypothesis is
$$
K_{m_0}: \;\limsup_{x\to\infty} x^{m_0}(1-F(x))<1.
$$
Notice that for $m_0=2$ the non-rejection $H_{m_0}$ implies that $F$ 
has infinite second moment. Although this test has a very good 
behavior even for small samples, it is only applicable to absolutely 
continuous distributions whereas our test does not need such a 
condition, which is essential for applications to samples produced by complex
simulations of random processes. 
Moreover,   Jure{\v{c}}ková and Picek's test applies to I.I.D. random variables,
whereas our methodology 
can be extended to some weak dependence cases (see Section~\ref{sec:num-exp-weak-dep} below). 
This potentially allows applications to interacting 
particles which propagate chaos. 
Finally, we refer to~\cite{markovich2008nonparametric} for other 
nonparametric techniques in the context of heavy tailed distributions.
\end{remark}

\section{Consistency of the statistic $\hS_n^m$: Proof of 
Prop.~\ref{prop:conv-in-P-estimator}} 
\label{sec:consistency}

The objective of this section is to establish the consistency 
proposition~\ref{prop:conv-in-P-estimator}.

To prove this proposition we need two preliminary results.

\subsection{Two consistency properties of $\hS_n$}
\begin{prop}\label{prop:Sn_cont}
{Recall that we have set $Z=\Psi(L)$.} 
For any $n\in\N$ it holds that
\begin{equation}\label{eq:conv-D-hSn}
{\hS_n^m} \Dconv{m\to\infty} \hS_n(Z).
\end{equation}
\end{prop}

\begin{proof}
For any positive integer $n$ recall the functional $\hS_n:D[0,1]\to\R$ defined in \eqref{eq:functional-hS}.
This functional is not continuous on $D[0,1]$ for the Skorokhod 
topology. Its discontinuity set is 
\begin{equation*}
\begin{aligned}
 J_n := &\{ z\in 
D[0,1]: \exists i=1,\ldots,n-1:\D z(i/n)>0 \} \\
&\cup  \{ z\in D[0,1]: 
\forall i=0,\ldots,n-1: z(i/n)=z((i-1)/n) \}.	
\end{aligned}
\end{equation*}
Therefore, as soon as 
$\P(Z\in J_n)=0$ we can apply the Continuous Mapping Theorem to get
\begin{equation*}
\hS_n(Z^m) \Dconv{m\to\infty} \hS_n(Z).
\end{equation*}

On the one hand, if the subjacent random variable $X$ belongs to 
$\DA(2)$, then $Z$ is a Brownian Bridge and therefore $\P(Z\in J_n)=0$. 

On the other hand, if $X$ belongs to $\DA(\alpha)$ with $o<\alpha<2$, 
then the limit process $Z$ is equal to $\Psi(L)$ where $L$ is a 
$\alpha$-stable process. Notice that $Z$ and $L$ hav the same jumps. 
Since the probability of $L$ having a jump at any fixed time is zero, 
it follows that $\P\left(Z \text{ has a jump at } i/n \right)=0$. In 
addition, we have
$$\P(Z\in J_n)\leq \sum_{i=1}^{n-1}\P\left(Z \text{ has a jump at } i/n 
\right) + \P\left(Z_{i/n}=0,\;i=1,\ldots,n-1\right) =0,$$
since since $Z_0=Z_1=0$, each summand in the first term of the 
right-hand side is zero, and the second term in the right-hand side is 
zero because $\law(Z_t)$ is absolutely continuous with respect to the 
Lebesgue measure (see Bertoin~\cite[p.218]{bertoin1998levy}).  

To summarize, in both cases we have $\P(Z\in J_n)=0$ and 
\eqref{eq:conv-D-hSn} holds true.
\end{proof}

\begin{prop}\label{prop:conv-in-P-Bivariation-Bridges} 
Let $\Psi$ be defined as in Lemma 
\ref{lem:continuity-of-Psi-Skokokhod-Topology}
and $\hS_n$ be defined as in~\eqref{def:S-hat-n}. Then,

\begin{enumerate}
\item If $(Y_t)_{0\leq t\leq 1}$ is a standard Brownian motion one has
$$\hS_n(\Psi(Y))\ASconv{n\to\infty}\frac{2}{\pi}.$$ 
\item If $(Y_t)_{0\leq t\leq 1}$ is an  $\alpha$-stable process 
starting from~$0$ one has
$$ \forall 1<\alpha\leq 2,~~\hS_n(\Psi(Y))\ASconv{n\to\infty}0. $$
For $\alpha\leq1$ the convergence holds in probability only. 
\end{enumerate} 
\end{prop} 

\begin{proof}
In view of \eqref{def:S-hat-n} we succesively consider $\hQ(\Psi(Y),n)$ and $\hB(\Psi(Y),n)$. 

A straightforward computation leads to
\begin{equation}\label{eq:hQPhi}
\hQ(\Psi(Y),n) =   \hQ(Y,n) -\frac{1}{n}Y_1^2.
\end{equation}
The second term in the right-hand side converges to zero a.s. We thus 
have to prove the a.s. convergence of~$\hQ(Y,n)$.

When $Y$ is a standard Brownian motion we have
$$   \hQ(Y,n) = \frac{1}{n}\sum_{i=1}^{n}\left( \sqrt{n}\D_i^nY \right)^2 \ASconv{n\to\infty} 1,$$ 
thanks to the Strong Law of Large Numbers, from which
\begin{equation}\label{eq:AsConvQuadraticVarBrownianBridge}
\hQ(\Psi(Y),n)   \ASconv{n\to\infty} 1 
\end{equation} 
and
$$\lim_{n\to\infty} \hS_n(\Psi(Y))=
\lim_{n\to\infty}\frac{\hB(\Psi(Y),n)}{\hQ(\Psi(Y),n)}
=\lim_{n\to\infty}\hB(\Psi(Y),n).$$
Consequently, the first statement of the proposition results from
the two following results which will be proven in the next subsection 
(see Lemmas~\ref{lem:bivariation-brownian-motion-alpha-stable} 
and~\ref{lem:convergence-of-bivariation-bridge}): For any Brownian 
motion~$Y$,
$$ \hB(Y,n)\ASconv{n\to\infty} \frac{2}{\pi} $$
and
$$\hB(Y,n)- \hB(\Psi(Y),n)\ASconv{n\to\infty}0.$$

Now, when $Y$ is an $\alpha$-stable process with $\alpha<2$ we have
$$ \hQ(Y,n) = \sum_{i=1}^{n}\left( \D_i^nY \right)^2 = 
\frac{1}{n^{2/\alpha}}\sum_{i=1}^{n}\left(n^{1/\alpha} \D_i^nY
\right)^2 $$ 
and the positive random variable $U_i =\left(n^{1/\alpha} \D_i^nY
\right)^2$ satisfies $\P(U_i>t)\sim t^{-\alpha/2}$.
Consequently, as $n$ goes to infinity, $\hQ(Y,n) $ converges in distribution 
 to a stable law $Q$ with index $\frac{\alpha}{2}<1$ (see Theorem~4.5.2 in~\cite[Sec.5, Chap.4]{whitt2002stochastic}), from which 
 $\P(Q=0)=0$. From \eqref{eq:hQPhi} it follows that~$\hQ(\Psi(Y),n)$ also converges
 in distribution to $Q$. We now use the two following results
which will be proven in the next subsection: If $Y$ is an 
$\alpha$-stable process starting from~0, then
$$ \hB(Y,n)\ASconv{n\to\infty} 0$$
and 
$$\lim_{n\to\infty}\hB(\Psi(Y),n)=\lim_{n\to\infty}\hB(Y,n), $$
where the limit is in probability.

Consequently, in view of Slutsky's lemma,
$$\hS_n(\Psi(Y))=
\frac{\hB(\Psi(Y),n)}{\hQ(\Psi(Y),n)}\Pconv{n\to\infty}0.$$ 
\end{proof}

In the proof of Proposition~\ref{prop:conv-in-P-Bivariation-Bridges} we 
have used the two following lemmas.

\subsection{Two key lemmas}

The first point of the first Lemma in this subsection is contained in 
Barndorff-Nielsen and Shephard~\cite[Thm.4.]{barndorff2006econometrics}.
We here give its easy proof for the sake of completeness. The second 
point requires less obvious arguments.

\begin{lem} \label{lem:bivariation-brownian-motion-alpha-stable} 
We have:
\begin{enumerate}
\item Let $(Y_t)_{0\leq t\leq 1}$ be a standard Brownian motion, then
$$\hB(Y,n)\ASconv{n\to\infty}\frac{2}{\pi}.$$
\item Let $(Y_t)_{0\leq t\leq 1}$ be an  $\alpha$-stable process 
starting from $0$. For $\alpha>1$, we have
$$\hB(Y,n)\ASconv{n\to\infty}0,$$
whereas for  $\alpha\leq1$ this last  convergence holds in probability 
only.
\end{enumerate}
\end{lem}

\begin{proof}

\begin{enumerate}
\item Let $(Y_t)_{0\leq t\leq 1}$ be a Brownian motion. Then,
\begin{align*}
 \sum_{i=1}^{n-1} |\D_i^n Y||\D_{i+1}^n Y|   & =    \sum_{i=1}^{\lfloor 
 n/2 \rfloor} |\D_{2i-1}Y | |\D_{2i} Y| + \sum_{i=1}^{\lfloor n/2 
 \rfloor-1} |\D_{2i}Y | |\D_{2i+1}Y |\\
  &\quad + 2\left(\frac{n}{2} - \lfloor \frac{n}{2} \rfloor 
  \right)|\D_{n-1}^n Y||\D_{n}^n Y|\\
  & =    \frac{1}{2}\frac{2}{n}\sum_{i=1}^{\lfloor n/2 \rfloor} 
  |\sqrt{n}\D_{2i-1}Y | |\sqrt{n}\D_{2i} Y| + \frac{1}{2}\frac{2}{n} 
  \sum_{i=1}^{\lfloor n/2 \rfloor-1} |\sqrt{n}\D_{2i}Y | 
  |\sqrt{n}\D_{2i+1}Y |\\
  &\quad + 2\left(\frac{n}{2} - \lfloor \frac{n}{2} \rfloor 
  \right)|\D_{n-1}^n Y||\D_{n}^n Y|.
\end{align*}
Notice that $\sqrt{n}\D_{i}Y\sim\mathcal{N}(0,1)$. The Law of the Large 
Numbers imples 
$$\frac{2}{n}\sum_{i=1}^{\lfloor n/2 \rfloor} |\sqrt{n}\D_{2i-1}Y | 
|\sqrt{n}\D_{2i} Y|  \ASconv{n\to\infty} \frac{2}{\pi}, $$
$$\frac{2}{n}\sum_{i=1}^{\lfloor n/2 \rfloor} |\sqrt{n}\D_{2i}Y | 
|\sqrt{n}\D_{2i+1} Y|  \ASconv{n\to\infty} \frac{2}{\pi}. $$
In addition, the continuity of the trajectories of the Brownian motion
implies that
 $$\left| 2\left(\frac{n}{2} - \lfloor \frac{n}{2} \rfloor 
 \right)|\D_{n-1}^n Y||\D_{n}^n Y| \right|\leq 2|\D_{n-1}^n Y||\D_{n}^n 
 Y|\ASconv{n\to \infty}0,$$
 Putting all together we get
 $$\hB(Y,n)\ASconv{n\to\infty}\frac{2}{\pi}.$$

 \item Let $(Y_t)_{0\leq t\leq 1}$ be an  $\alpha$-stable process 
 starting from $0$. In this case we have
\begin{align*}
\hB(Y,n)  &=      \frac{1}{2}\frac{2}{n^{2/\alpha}}\sum_{i=1}^{\lfloor 
n/2 \rfloor} |n^{1/\alpha}\D_{2i-1}Y | |n^{1/\alpha}\D_{2i} Y| + 
\frac{1}{2}\frac{2}{n^{2/\alpha}} \sum_{i=1}^{\lfloor n/2 \rfloor-1} 
|n^{1/\alpha}\D_{2i}Y | |n^{1/\alpha}\D_{2i+1}Y |\\
  &\quad + 2\left(\frac{n}{2} - \lfloor \frac{n}{2} \rfloor 
  \right)|\D_{n-1}^n Y||\D_{n}^n Y|,    
\end{align*}
and for any $i=1,\ldots,n$, $n^{1/\alpha}\D_{i}Y$ are independent 
$\alpha$-stable random variables. We proceed by cases:
\begin{itemize}
\item $\alpha\in(1,2):$ Notice that
$$ \frac{2}{n^{2/\alpha}}\sum_{i=1}^{\lfloor n/2 \rfloor} 
|n^{1/\alpha}\D_{2i-1}Y || n^{1/\alpha}\D_{2i} Y|  
= \frac{1}{n^{2/\alpha-1}}\left(\frac{2}{n}
\sum_{i=1}^{\lfloor n/2 \rfloor} |n^{1/\alpha}
\D_{2i-1}Y | |n^{1/\alpha}\D_{2i} Y| \right). $$
The Law of Large Numbers implies that the term inside the parenthesis in
the right-hand side a.s. converges to a finite quantity. Therefore the 
whole right-hand side a.s. tends to~0. Similarly,
$$\frac{2}{n^{2/\alpha}} \sum_{i=1}^{\lfloor n/2 \rfloor-1}
|n^{1/\alpha}\D_{2i}Y | |n^{1/\alpha}\D_{2i+1}Y |
\ASconv{n\to\infty}0.$$
We finally observe that 
$$\limsup_{n\to\infty } \left|2\left(\frac{n}{2} - \lfloor \frac{n}{2}
\rfloor \right)|\D_{n-1}^n Y||\D_{n}^n Y|\right| 
\leq 2 \lim_{n\to\infty}\left|Y_{1} - Y_{1-\frac{1}{n}} \right| \left|Y_{1-\frac{1}{n}} - Y_{1-\frac{2}{n}} \right| =0\;\;\text{a.s.} $$
since the paths of $Y$ have left limits. Therefore, for $\alpha\in(1,2)$ we have
$$\hB(Y,n)\ASconv{n\to\infty}0.$$

\item $\alpha\in(0,1):$ In view of 
Embrechts and Goldie~\cite[Cor. of Thm.3]{embrechts1980closure} we have 
that  
$$|n^{1/\alpha}\D_{2i-1}Y| | n^{1/\alpha}\D_{2i} Y|\in 
~\DA(\alpha).$$ Therefore, there exists a slowly varying 
function~$\ell_0$ such that
$$  \frac{1}{n^{1/\alpha}\ell_0(n)}
\sum_{i=1}^{\lfloor n/2\rfloor} |n^{1/\alpha}\D_{2i-1}Y
||n^{1/\alpha}\D_{2i} Y|, $$
converges in distribution to an $\alpha$-stable random variable. Hence
$$ \frac{2}{n^{2/\alpha}}\sum_{i=1}^{\lfloor n/2 \rfloor} 
|n^{1/\alpha}\D_{2i-1}Y || n^{1/\alpha}\D_{2i} Y| =
\frac{\ell_0(n)}{n^{1/\alpha}}\frac{1}{n^{1/\alpha}\ell_0(n)}
\sum_{i=1}^{\lfloor n/2 \rfloor} |n^{1/\alpha}\D_{2i-1}Y| |
n^{1/\alpha}\D_{2i} Y|$$
tends to zero in probability, because is the product of a sequence of
random variables converging in distribution and a real sequence
converging to zero. For the other terms we can proceed as before.
We thus are in a position to conclude that
$$\hB(Y,n)\Pconv{n\to\infty}0.$$ 
\item $\alpha=1:$ We need to
study normalized sums of the type 
$$\frac{1}{2}\frac{2}{n^2}\sum_{i=1}^{\lfloor n/2 \rfloor} 
|n\D_{2i-1}Y || n\D_{2i} Y|.$$ 
Denote by $A$ such a quantity. Notice that for each $i$,  $n\D_{i}Y$ 
has a  $1$-stable law. By again using Embrechts and 
Goldie~\cite[Cor. of Thm.3]{embrechts1980closure}
we get that $ |n\D_{2i-1}Y | n\D_{2i} Y|$ belongs to~$\DA(1)$. 
Hence, we need to introduce centering factors to obtain the convergence 
of~$A$. 

In view of Whitt~\cite[Thm.4.5.1]{whitt2002stochastic} the centered and 
normalized sums
\begin{align*}
 \frac{1}{(n/2)\ell_0(n/2)}  &\Bigg(\sum_{i=1}^{\lfloor n/2 \rfloor} 
 |n\D_{2i-1}Y || n\D_{2i} Y|\\
 &\quad - (n/2)^2\ell_0(n/2) 
 \E\left[\sin\left(\frac{|n\D_{2i-1}Y || n\D_{2i} 
 Y|}{(n/2)\ell_0(n/2)}\right)\right]  \Bigg) 
 \end{align*}
converge in distribution to a $1$-stable random variable. Then

\begin{align*}
A &=      \frac{\ell_0(n/2)}{2n}\Big( 
\text{Converging in Distribution Sequence} \Big)  \\
& \quad+\ell_0(n/2) \E\left[\sin\left(\frac{|n\D_{1}Y | |n\D_{2}
 Y|}{(n/2)\ell_0(n/2)}\right)\right] ,
\end{align*}
The first term tends to $0$ in probability because is the product of a
deterministic sequence converging to $0$ and a weakly convergent
sequence of random variables. As for the second one, there
exist $1$-stable random variables $S_1$ and $S_2$ such that
\begin{align*} 
\ell_0(n/2) \E\left[\sin\left(\frac{|n\D_{2i-1}Y ||
n\D_{2i} Y|}{(n/2)\ell_0(n/2)}\right)\right]  &=    \ell_0(n/2)
\E\left[\sin\left(\frac{|S_1S_2|}{(n/2)\ell_0(n/2)}\right)\right] \\ 
&\leq 
\E\left[\sqrt{|S_1|}\right]^2\frac{\sqrt{\ell_0(n/2)}}{\sqrt{n/2}},
\end{align*} 
where in the last inequality we have used $|\sin(x)|\leq
\sqrt{|x|}$. Hence the right-hand side converges to~$0$. We thus
conclude that $$\hB(Y,n)\Pconv{n\to\infty}0.$$ 
\end{itemize}

\end{enumerate} 
\end{proof}

The second lemma in this subsection shows that $\hB(Y,n)$ and 
$\hB(\Psi(Y),n)$ have a similar asymptotic behaviour.

\begin{lem}\label{lem:convergence-of-bivariation-bridge}
If $(Y_t)_{0\leq t\leq 1}$ is, either a Brownian motion or an
$\alpha$-stable process starting from $0$ with $\alpha>1$, then
$$\hB(Y,n)- \hB(\Psi(Y),n)\ASconv{n\to\infty}0.$$ 
Moreover, if  $(Y_t)_{0\leq t\leq 1}$ is an
$\alpha$-stable process starting from $0$ with   $\alpha\leq1$ this 
last  convergence holds in probability only.
\end{lem}

\begin{proof} [Proof of Lemma 
\ref{lem:convergence-of-bivariation-bridge} ]
Since for all $a,b\in\R$ $||a|-|b||\leq |a-b|$, we have
\begin{align*} 
\left| \sum_{i=1}^{n-1} |\D_i^n \Psi(Y)||\D_{i+1}^n
\Psi(Y)| -  \sum_{i=1}^{n-1} |\D_i^n Y||\D_{i+1}^nY|\right|     & \leq 
\sum_{i=1}^{n-1} \left| \D_i^n \Psi(Y)\D_{i+1}^n \Psi(Y) - \D_i^n
Y\D_{i+1}^nY\right|,   
\end{align*} 
and then
\begin{equation}
\label{eq:firts-bound-lemma-convergence-of-bivariation-bridges}
\begin{aligned} 
\left| \sum_{i=1}^{n-1} |\D_i^n \Psi(Y)||\D_{i+1}^n
\Psi(Y)| -  \sum_{i=1}^{n-1} |\D_i^n Y||\D_{i+1}^nY|\right|   &\leq 
\sum_{i=1}^{n-1} \left|   \frac{Y_1}{n} \D_i^n Y  \right|+
\sum_{i=1}^{n-1} \left|\frac{Y_1}{n} \D_{i+1}^nY  \right|\\
&\quad+
\sum_{i=1}^{n-1} \left|  \frac{1}{n^2}Y_1^2 \right|. 
\end{aligned} 
\end{equation}

We now analyze each term in the right-hand side of 
\eqref{eq:firts-bound-lemma-convergence-of-bivariation-bridges}. We 
start with the third one: 
$$ \sum_{i=1}^{n-1} \left| \frac{1}{n^2}Y_1^2 \right| 
=   \frac{(n-1)}{n^2}Y_1^2 \ASconv{n\to\infty}0. $$
As for the second term in the right-hand side of  
\eqref{eq:firts-bound-lemma-convergence-of-bivariation-bridges}, we have
$$ \sum_{i=1}^{n-1} \left|   \frac{Y_1}{n} \D_i^n Y  \right| 
= \frac{|Y_1|}{n^{1/\alpha}}\frac{1}{n}  \sum_{i=1}^{n-1} 
\left| n^{1/\alpha}   \D_i^n Y  \right|. $$
Notice that,  for $i=1,\ldots,n-1$, $|n^{1/\alpha}   \D_i^n Y|$ are 
i.i.d. random variables with stable law
of index $\alpha$.

When~$Y$ is a Brownian motion or when $\alpha>1$, there
exists $y\in \R$ such that 
$\E\left[|n^{1/\alpha} \D_i^n Y|\right]=y$ for any $n\geq1$. Then, 
the Law of Large Numbers implies
$$\frac{1}{n}  \sum_{i=1}^{n-1} 
\left| n^{1/\alpha}   \D_i^n Y  \right|\ASconv{n\to\infty} y, $$
from which
$$ \sum_{i=1}^{n-1} \left|   \frac{Y_1}{n} \D_i^n Y  \right| 
= \frac{|Y_1|}{n^{1/\alpha}}\frac{1}{n}  \sum_{i=1}^{n-1} 
\left| n^{1/\alpha}   \D_i^n Y  \right|\ASconv{n\to\infty} 0.
$$

When $\alpha<1$ we use
Embrechts et al.~\cite[Thm.2.2.15]{embrechts2013modelling} to get that
there exists  a slowly varying function $\ell_0$ such that 
$$ \frac{1}{\ell_0(n)}  
\sum_{i=1}^{n-1} \left| \D_i^n Y  \right|$$ 
converges in distribution as $n\to\infty$ to a stable random variable. 
Therefore,
\begin{align*} \sum_{i=1}^{n-1} \left|   \frac{Y_1}{n} \D_i^n Y  \right|
&= \frac{|Y_1\ell_0(n)|}{n}\frac{1}{\ell_0(n)} 
\sum_{i=1}^{n-1} \left|\D_i^n Y  \right|, 
\end{align*}
is the product of a random variable which a.s. converges to~0 when~$n$ 
goes to infinity and of a random variable which 
converges in distribution, and thus converges to~0 in probability.

Finally, if $\alpha=1$, we have that $n\D_i^n Y\in\DA(1)$. Therefore,
in view of Whitt~\cite[Thm.4.5.1]{whitt2002stochastic}, 
there exists a slowly varying function $\ell_0$ such that 
$$\frac{1}{n\ell_0(n)} \left( \sum_{i=1}^{n}
|n\D_i^n Y|   - n^2\ell_0(n)\E\left[\sin\left(\frac{n\D_1^n
Y}{n\ell_0(n)}\right)\right]\right) $$
converges in distribution to a stable distribution. We conclude that
\begin{align*}
\sum_{i=1}^{n-1} \left|   \frac{Y_1}{n} \D_i^n Y  \right| &=  
\frac{\ell_0(n)|Y_1|}{n}\frac{1}{n\ell_0(n)} \left( \sum_{i=1}^{n}
|n\D_i^n Y|   - n^2\ell_0(n)\E\left[\sin\left(\frac{n\D_1^n
Y}{n\ell_0(n)}\right)\right]\right)   \\ 
& \quad
+\ell_0(n)\E\left[\sin\left(\frac{n\D_1^n Y}{n\ell_0(n)}\right)\right]
\end{align*} 
 converges to zero in probability since the first term in the 
 right-hand side is the product of a random variable which converges
 to zero a.s. when $n$ goes to infinity and a a random variables whcih 
 converges in distribution. As for the second term in the right-hand 
 side, one 
 can check that it goes to zero by proceeding 
 as in the end of the proof of the previous proposition. To summarize, 
 we have
$$ \sum_{i=1}^{n-1} \left|  
\frac{Y_1}{n} \D_i^n Y  \right|\Pconv{n\to\infty}0.$$ 
We thus have obtained that 
$$ \sum_{i=1}^{n-1} \left|   \frac{Y_1}{n} \D_{i+1}^n Y 
\right|\Pconv{n\to\infty}0 $$ 
and 
$$\left| \sum_{i=1}^{n-1} 
|\D_i^n\Psi(Y)||\D_{i+1}^n \Psi(Y)| -  \sum_{i=1}^{n-1} 
|\D_i^n Y|~|\D_{i+1}^nY|\right| \Pconv{n\to\infty}0. $$ 
\end{proof}

\subsection{End of the proof of 
Proposition~\ref{prop:conv-in-P-estimator}}
We now are in a position to prove the
proposition~\ref{prop:conv-in-P-estimator}.

We start with proving that for any bounded and continuous 
function~ $\varphi:\R\to\R$ we have 
\begin{equation} \label{limit-psi(kappa)}
\lim_{n\to\infty}\lim_{m\to\infty}
\E\left[\varphi({\hS_n^m})\right] = \varphi(\kappa).
\end{equation}
In view of Proposition \ref{prop:Sn_cont} we have: 
$$\lim_{m\to\infty}\left|\E\left[\varphi({\hS_n^m})\right] -
\E\left[\varphi(\hS_n(Z))\right] \right|=0. $$ 
In addition, in view of
Proposition~\ref{prop:conv-in-P-Bivariation-Bridges} we have
$$\lim_{n\to\infty}\left|\E\left[\varphi(\hS_n(Z))\right] - \varphi(\kappa)
\right|=0. $$ 
Therefore, the desired limit in~\eqref{limit-psi(kappa)} is obtained
by letting $m$ and then $n$ tend to infinity in the right-hand side of
$$\left|\E\left[\varphi({\hS_n^m})\right] -
\varphi(\kappa) \right| \leq      
\left|\E\left[\varphi({\hS_n^m})\right] -
\E\left[\varphi(\hS_n(Z))\right] \right| +  
\left|\E\left[\varphi(\hS_n(Z))\right] - \varphi(\kappa) \right|. $$ 

Now, fix $\epsilon>0$ and define the function $\psi$ as follows:
{$\psi$ takes the value~$0$ outside the closed interval
$[\kappa-\epsilon,\kappa+\epsilon]$,
growths linearly from 0~to~1 on $[\kappa-\epsilon,\kappa]$
and decreases linearly from 1~to~0 on $[\kappa,\kappa+\epsilon]$.}
From
$$\P\left(|\hS_n^m-\kappa|\leq \epsilon\right)
=\E\left[\ind{[\kappa-\epsilon,\kappa+\epsilon]}(\hS_n^m)\right]\geq
\E\left[\varphi(\hS_n^m)\right], $$ 
it results that
$$ \lim_{n\to\infty}\lim_{m\to\infty}\P\left(|\hS_n^m-\kappa|
\leq\epsilon\right) ]\geq 
\lim_{n\to\infty}\lim_{m\to\infty}\E\left[\varphi(\hS_n^m)\right] =
\varphi(\kappa)=1. $$

\section{A CLT for $\hS_n^m$: Proof of 
Prop.~\ref{prop:central-limit-theorem-statistic}}
\label{sec:CLT}

In this section we prove a Central Limit theorem for 
$\hS_n^m$ defined as 
in~Proposition~\ref{prop:central-limit-theorem-statistic}.

We start with proving the CLT for the pair $(\hB(Y,n),\hQ(Y,n))$, 
with $(Y_t)_{0\leq t\leq 1}$ being  a Brownian Motion.

\subsection{Three CLT for $(\hB(Y,n),\hQ(Y,n))$}

The first result in this subsection concerns the case where~$Y$ is a 
Brownian motion.  It is a particular case of 
Barndorff-Nielsen and Shephard~\cite[Thm.3]{barndorff2006econometrics}.
For the sake of completeness we here provide a straightforward proof
adapted to the specific Brownian case.

\begin{lem}\label{lem:bivariationBrownian} 
Let $(Y_t)_{0\leq t\leq 1}$be a Brownian motion. One then has 
$$ \left(\sqrt{n}\left[\hB(Y,n)-\frac{2}{\pi}\right] ,
\sqrt{n}\left[\hQ(Y,n)-1\right] \right)\Dconv{n\to\infty}
\mathcal{N}_2\left(0,\Sigma\right),$$ 
where
\begin{equation}\label{eq:CovarianceLimit} 
\Sigma= \left(
\begin{array}{cc} 
1 +\frac{4}{\pi}- \frac{12}{\pi^2} &  \frac{4}{\pi} \\
\frac{4}{\pi} &      2 
\end{array} 
\right). 
\end{equation}
\end{lem}

\begin{proof} 
Notice that
\begin{align*}
    \Bigg(\sqrt{n}&\left[\hB(Y,n)-\frac{2}{\pi}\right] ,
    \sqrt{n}\left[\hQ(Y,n)-1\right] \Bigg)\\
    &= \frac{1}{\sqrt{n}}
    \Bigg(\sum_{i=1}^{n-1} |\sqrt{n}\D_i^n Y||\sqrt{n}\D_{i+1}^n Y| -
    \frac{2}{\pi} , \sum_{i=1}^{n-} |\sqrt{n}\D_i^n Y|^2 -1\Bigg).
\end{align*}

The right-hand side has the same law as 
\begin{align*}
\frac{1}{\sqrt{n}}\left(  \sum_{i=1}^{n-1} \left[|G_i ||G_{i+1} | - 
\frac{2}{\pi}\right], \sum_{i=1}^{n} \left[G_i^2 -  1\right] \right),
\end{align*} 
where $G_i$ is a standard Normal random variable.

The vectors $(|G_i ||G_{i+1} | -  {2}/{\pi},G_i^2 -  1 )$ are 
identically distributed, centered, and such that   $(|G_i
||G_{i+1} | -  {2}/{\pi},G_i^2 -  1 ) $ is independent of  $(|G_j
||G_{j+1} | -  {2}/{\pi},G_j^2 -  1 ) $ if $|j-i|>1$. Hence, we can
apply the Central Limit Theorem for $1$-dependent sequences of vectors
(see Hoeffding and Robbins~\cite[Thm.3]{hoeffding1948}). The
asymptotic covariance matrix is given by 
\begin{equation*} 
\Sigma_{11}:=
    		\E\left[\left(|G_{i+1} ||G_{i+2} | -  \frac{2}{\pi} \right)^2 \right] + 
    		2\E\left[\left(|G_i ||G_{i+1} | -  \frac{2}{\pi} \right)\left(|G_{i+1}
   		 ||G_{i+2} | -  \frac{2}{\pi} \right) \right]     
    =  1 +\frac{4}{\pi}- \frac{12}{\pi^2}, 
\end{equation*}
\begin{equation*} 
    \Sigma_{12}
    := \E\left[\left(|G_{1} ||G_{2} | - 
        \frac{2}{\pi} \right)\left(G_1^2 -  1\right) \right] + 
        \E\left[\left(|G_1 ||G_{2} | -  \frac{2}{\pi} \right) \left(G_2^2 - 
        1\right)  \right]   +  \E\left[\left(|G_2 ||G_{3} | -  \frac{2}{\pi}
        \right) \left(G_1^2 -  1\right)  \right]
    =\frac{4}{\pi},
 \end{equation*} 
and
\begin{equation*} 
\Sigma_{22} := \E\left[\left(G_1^2 - 1\right)^2 \right] 
+ 2\E\left[\left(G_1^2 -  1\right) \left(G_2^2 - 1\right)  \right]  
=2. 
\end{equation*} 
\end{proof}

In the next Lemma we extend the previous result to Brownian Bridges.
 
\begin{lem}\label{lem:CLTBrownianBridge}
Let $(Y_t)_{0\leq t\leq 1}$ be a Brownian motion. Then
\begin{equation} \label{L2-estimate-B-hat}
\E\left[\left(\sqrt{n}\left[\hB(\Psi(Y),n)-\hB(Y,n)\right]\right)^2
\right]+\E\left[\left(\sqrt{n}\left[\hQ(\Psi(Y),n)-\hQ(Y,n)\right]
\right)^2\right]\leq\frac{C}{n}.
\end{equation}
In addition, one has
\begin{equation}\label{eq:D-convergence-B-hat-Q-hat}
\left(\sqrt{n}\left[\hB(\Psi(Y),n)-\frac{2}{\pi}\right] ,  
\sqrt{n}\left[\hQ(\Psi(Y),n)-1\right]  \right)\Dconv{n\to\infty} 
\mathcal{N}_2\left(0,\Sigma\right),
\end{equation}
where $\Sigma$ is defined as in~\eqref{eq:CovarianceLimit}. 
\end{lem}

\begin{proof} A straightforward computation shows that
$$ \hQ(\Psi(Y),n)= \hQ(Y,n)-\frac{1}{n}Y_1^2. $$
Therefore, we have
$$\E\left[\left(\sqrt{n}\left[\hQ(\Psi(Y),n)
-\hQ(Y,n)\right]\right)^2\right]
 = \E\left[\frac{1}{n}Y_1^4\right]=\frac{3}{n}.$$ 
Now consider a family of independent random variables 
$G_i\sim\mathcal{N}(0,1)$ and set
$\bar{G}^n:=n^{-1}\sum_{i=1}^{n}G_i$. It holds that
\begin{align*} 
\E\left[\left(\sqrt{n}\left[\hB(\Psi(Y),n)-\hB(Y,n)\right]\right)^2\right]
 &= n \sum_{i=1}^{n-1}\sum_{k=1}^{n-1}\E\bigg[ \left(|\D_i^n 
\Psi(Y)||\D_{i+1}^n \Psi(Y)| - |\D_i^n Y||\D_{i+1}^n Y|\right)\\ 
&\qquad\qquad\times\left( |\D_k^n \Psi(Y)||\D_{k+1}^n \Psi(Y)| - 
|\D_k^n Y||\D_{k+1}^n Y|\right)  \bigg]  \\ &= \frac{1}{n} 
\sum_{i=1}^{n-1}\sum_{k=1}^{n-1}\E\bigg[\left(|G_i - \bar{G}^n||G_{i+1} 
- \bar{G}^n| - |G_{i} ||G_{i+1} |\right)\\ 
&\qquad\qquad\times\left(|G_{k} - \bar{G}^n||G_{k+1} - \bar{G}^n| - 
|G_{k} ||G_{k+1} |\right)\bigg]  \\ & =: S_1+S_2, 
\end{align*} 
where we have split the sum according to $k=i,i+1$ in $S_1$ and 
$k\notin\{i,i+1\}$ in $S_2$.

As for $S_1$, notice that $$ |G_i - \bar{G}^n||G_{i+1} - \bar{G}^n| -
|G_{i} ||G_{i+1} | \leq |\bar{G}^n||\bar{G}^n- G_i - G_{i+1}|,$$ and
therefore 
\begin{align*} 
S_1 & =  \frac{1}{n}
\sum_{i=1}^{n-2}\sum_{k=i}^{i+1}\E\left[ \left(|G_i - \bar{G}^n||G_{i+1}
- \bar{G}^n| - |G_{i} ||G_{i+1} |\right)\left(|G_{k} -
\bar{G}^n||G_{k+1} - \bar{G}^n| - |G_{k} ||G_{k+1} |\right) \right]  \\
&\leq  \frac{1}{n} \sum_{i=1}^{n-2}\sum_{k=i}^{i+1}\sqrt{\E\left[
\left(|G_i - \bar{G}^n||G_{i+1} - \bar{G}^n| - |G_{i} ||G_{i+1}
|\right)^2\right]}\\ &\quad\quad\times\sqrt{\E\left[ \left(|G_{k} -
\bar{G}^n||G_{k+1} - \bar{G}^n| - |G_{k} ||G_{k+1} |\right)^2 \right] }
\\ & \leq   \frac{1}{n} \sum_{i=1}^{n-2}\sum_{k=i}^{i+1}\sqrt{\E\left[
|\bar{G}^n|^2|\bar{G}^n- G_i - G_{i+1}|^2\right]\E\left[ 
|\bar{G}^n|^2|\bar{G}^n- G_{k} - G_{k+1}|^2 \right] }   \\ 
& \leq   
\frac{1}{n} \sum_{i=1}^{n-2}\sum_{k=i}^{i+1}\sqrt[4]{\E\left[
|\bar{G}^n|^4\right]\E\left[|\bar{G}^n- G_i - G_{i+1}|^4 \right]\E\left[
|\bar{G}^n|^4\right]\E\left[|\bar{G}^n- G_{k} - G_{k+1}|^4 \right]} \\
&\leq \frac{C}{n}. 
\end{align*} 
As for $S_2$, we use the following inequality which we admit for a 
while (see Subsection~\ref{proof-key-inequality} for its proof):
\begin{equation}\label{eq:inequality-CLT-for-brownian-bridge}
\E\bigg[\left(|G_i - \bar{G}^n||G_{i+1} - \bar{G}^n| - |G_{i} ||G_{i+1}
|\right)\left(|G_{k} - \bar{G}^n||G_{k+1} - \bar{G}^n| - |G_{k}
||G_{k+1} |\right)\bigg] \leq \frac{C}{n^2}.
\end{equation} 
This ends the proof of~\eqref{L2-estimate-B-hat}.

From ~\eqref{L2-estimate-B-hat} it follows that 
$(\sqrt{n}[\hB(\Psi(Y),n)-\hB(Y,n)] 
,\sqrt{n}[\hQ(\Psi(Y),n))-\hQ(Y,n))]$ converges in probability to zero. 
Therefore, Slutsky's Theorem implies that
\begin{equation*}
\begin{aligned}
\left(\sqrt{n}\left[\hB(\Psi(Y),n)-\frac{2}{\pi}\right]\right. &,  \left.
\sqrt{n}\left[\hQ(\Psi(Y),n)-1\right]  \right)\\
 &=
\left(\sqrt{n}\left[\hB(\Psi(Y),n)-\hB(Y,n)\right] ,
\sqrt{n}\left[\hQ(\Psi(Y),n))-\hQ(Y,n))\right]\right)  \\		\\
&\quad +  \left(\sqrt{n}\left[\hB(Y,n)-\frac{2}{\pi}\right] ,  
\sqrt{n}\left[\hQ(Y,n)-1\right]  \right) 		
\end{aligned}
\end{equation*}
converges in distribution to $\mathcal{N}_2(0,\Sigma)$.

\end{proof}
We thus now in a position to obtain the following result.
\begin{prop}
\label{prop:asymptotic-distribution-statistic-for-brownian-bridge}
Let $(Y_t)_{0\leq t\leq 1}$ be a Brownian motion. One then has:
$$\sqrt{n}\left[\hS_n(\Psi(Y))-\frac{2}{\pi}\right]\Dconv{n\to\infty}
\mathcal{N}\left(0, 1 +\frac{4}{\pi}- \frac{20}{\pi^2} \right).$$
\end{prop} 

\begin{proof}
Notice that
\begin{align*}
\sqrt{n}\left[\frac{\hB(\Psi(Y),n)}{\hQ(\Psi(Y),n)}-\frac{2}{\pi}\right]
&= \frac{\sqrt{n}}{\hQ(\Psi(Y),n)}\left[\hB(\Psi(Y),n)
-\frac{2}{\pi}\hQ(\Psi(Y),n)\right] \\ 
& =   \frac{1}{\hQ(\Psi(Y),n)}\left(\sqrt{n}
\left[\hB(\Psi(Y),n)-\frac{2}{\pi}\right]
-\frac{2}{\pi}\sqrt{n}\left[\hQ(\Psi(Y),n)-1\right]\right). 
\end{align*} 
In view of Lemma~\ref{lem:CLTBrownianBridge} 
the term in the parenthesis in the right-hand side converges to a 
centered Gaussian random variable with variance 
$$\Sigma_{11} - 2\frac{2}{\pi}\Sigma_{12} +\frac{4}{\pi^2} \Sigma_{22} 
= 1 +\frac{4}{\pi}- \frac{12}{\pi^2} 
- 2\frac{2}{\pi}\times\frac{4}{\pi} + \frac{4}{\pi^2}\times2 
= 1+\frac{4}{\pi}- \frac{20}{\pi^2}. $$

In addition, $\hQ(\Psi(Y),n)\ASconv{n\to\infty}1$ in view of 
\eqref{eq:AsConvQuadraticVarBrownianBridge}. Therefore, Slutky's 
Theorem allows us to conclude.
\end{proof}

\subsection{End of the proof 
of~Prop.\ref{prop:central-limit-theorem-statistic}}

We now are in a position to end the proof 
of~Prop.\ref{prop:central-limit-theorem-statistic}. Observe that for $Y$ a Brownian motion
\begin{align*}
&\left|\E\left[\varphi\left(\frac{\sqrt{n}}{\sigma_\pi}\left({\hS_n^m}
-\frac{2}{\pi}\right)\right)\right]
  - \E\left[\varphi(\mathcal{N})\right] \right| \\ 
& \quad\quad\quad\quad   \leq   \left|  
\E\left[\varphi\left(\frac{\sqrt{n}}{\sigma_\pi}
\left({\hS_n^m}-\frac{2}{\pi}\right)\right)\right]
- \E\left[\varphi\left(\frac{\sqrt{n}}{\sigma_\pi}\left(\hS_n(\Psi(Y))
-\frac{2}{\pi}\right)\right)\right]\right|  \\ 
&\quad\quad\quad\quad\quad + 
\left|\E\left[\varphi\left(\frac{\sqrt{n}}{\sigma_\pi}\left(\hS_n(\Psi(Y))
-\frac{2}{\pi}\right)\right)\right]
- \E\left[\varphi(\mathcal{N})\right] \right|.
\end{align*} 
Thanks to the Mapping Theorem 
and Proposition \ref{cor:convergenceToBridge}, we get that the first 
term in the right-hand side goes to zero when $m$ tends to infinity.  Thanks to 
Proposition~\ref{prop:asymptotic-distribution-statistic-for-brownian-bridge}, we get that the second term in the 
right-hand side also tends to zero when $n$ tends to 
infinity.

\subsection{A corollary to Tanaka's formula}
Before proving 
Inequality~\eqref{eq:inequality-CLT-for-brownian-bridge}
we need to deduce a key estimate from Tanaka's formula.

\begin{prop} \label{prop:ItoTanakaF} 
Let $(W_s,s\geq0)$ a standard
Brownian motion, $G\in\R^{2d}$ a random vector with density 
$f_G$ with respect
to the Lebesgue measure and independent of $W$ and
$f:\R\times\R^{2d}\to\R$. Let us denote
$\zi=(z_1,\ldots,z_{i-1},z_{i+1},\ldots,z_{2d})\in\R^{2d-1}$,
$\zi(x)=(z_1,\ldots,z_{i-1},x,z_{i+1},\ldots,z_{2d})\in\R^{2d}$ and
$$f_z:\R\to\R,\;\; x\to f_z(x):=f(x,z).$$ 

Suppose:
\begin{enumerate} \item For
any $z\in\R^{d}$, $f_z$ is regular enough to apply the It\^o-Tanaka
formula
\begin{equation}\label{eq:ItoTanaka-fz}
f_z(W_t)=f_z(0)+\int_{0}^{t}D_{-}f_z(W_s)~dW_s
+\frac{1}{2}\int_{-\infty}^{\infty}L_t^y(W)~D^2f_z(dy),
\end{equation} 
 where $D_{-}f_z$ is the left derivative of $f_z$, 
$D^2f_z$ is the second derivative of $f$ in the distributional sense 
and  $L_t^y(W)$ is the Local time of $W$. 
\item For any $z\in\R^{d}$, 
$f_z(0)=0$. 
\item For any $z\in\R^{d}$, $D_{-}f_z$ has polynomial 
growth such that the stochastic integral in \eqref{eq:ItoTanaka-fz} is 
a martingale. 
\item For any $z\in\R^{d}$, the measure $D^2f_z$ can be decomposed  
as 
\begin{equation}\label{eq:decomposition-second-derivative-f}
D^2f_z(dy) = f_z''(y)~dy + \sum_{i=1}^{2d}\gamma_z^i 
~\delta_{z_i}(dy), 
\end{equation}
where $ \delta_{a}$ is de Dirac distribution at point~$a\in\R$. 
\item $\E\left[f_G''(W_t)\right]$ is bounded  in 
a neighborhood of $t=0$. 
\item For any $i=1,\ldots,2d$,  
$\E\left[h_i(W_t)\right]$ is bounded in a neighborhood of $t=0$, where 
$$h_i(x)=\left(\int_{\R^{2d-1}}\gamma_{\zi(x)}^i~f_G(\zi|x)~d\zi\right) 
~f_{G_i}(x), $$ 
where the $\gamma_{\cdot}^i$'s are given by 
\eqref{eq:decomposition-second-derivative-f}, $f_G(\zi|x)$ is the 
conditional density of the vector $G$ knowing that~$G_i=x$
and $f_{G_i}$ is the Gaussian density of~$G_i$.
\end{enumerate} 

Then one has
$\E\left[f(W_t,G)\right]=O(t)$
for any $t$ in a neighborhood of $0^+$. 

Moreover, assume now that 
conditions 5 and 6 are respectively replaced by 
\begin{enumerate} 
\item[5'.] The 
function $x\to h_0(x):= \E\left[f_G''(x)\right]$ is null at zero, twice 
differentiable with first and second derivative that have at most 
polynomial growth. \item[6'.] For any $i=1,\ldots,2d$, $h_i(0)=0$, 
$h_i\in\mathcal{C}^2$, and $h_i'$ and $h_i''$ have at most polynomial 
growth.
\end{enumerate} 
Then, it holds that
$\E\left[f(W_t,G)\right]=O(t^2)$ for any $t$ in a neighborhood of $0^+$.  
\end{prop} 

\begin{proof}
From Conditions 1, 2 and 3 if follows that
$$\E\left[f_z(W_t)\right]=\frac{1}{2}\E\left[\int_{-\infty}^{\infty}
D^2f_z(dy)L_t^y(W)\right]. $$ 
From condition 4 it then results that
$$ \E\left[f_z(W_t)\right]=\frac{1}{2}\E\left[\int_{-\infty}^{\infty}
f_z''(y)L_t^y(W)dy\right]+\frac{1}{2}\E\left[\sum_{i=1}^{2d}
\gamma_z^iL_t^{z_i}(W)\right]. $$ 
By using the Occupation Time Formula we thus get
$$ \E\left[f_z(W_t)\right]
=\frac{1}{2}\int_{0}^{t}\E\left[f_z''(W_s)\right]ds
+\frac{1}{2}\E\left[\sum_{i=1}^{2d}\gamma_z^iL_t^{z_i}(W)\right]. $$ 
In addition, as $G$ and $W$ are independent we have
\begin{align*} 
\E\left[f(W_t,G)\right] 
&=\E\left[\E\left[f(W_t,z)\right]|_{z=G}\right]        \\ 
&=\E\left[\E\left[f_z(W_t)\right]|_{z=G}\right]        \\ 
& = \frac{1}{2}\int_{0}^{t}\E\left[f_G''(W_s)\right]ds
+\frac{1}{2}\E\left[\sum_{i=1}^{2d}\gamma_G^iL_t^{G_i}(W)\right]. 
\end{align*} 
Now observe that
\begin{align*} 
\E\left[\gamma_G^iL_t^{G_i}(W)\right] &= 
\E\int_{\R^{2d}}\gamma_z^iL_t^{z_i}(W)f_G(z)dz     \\ 
& = \E \int_{-\infty}^{\infty}\int_{\R^{2d-1}}\gamma_{\zi(z_i)}^i
L_t^{z_i}(W)f_G(\zi|z_i)f_{G_i}(z_i)d\zi dz_i     \\ 
& =  \E \int_{-\infty}^{\infty}\left(\int_{\R^{2d-1}}\gamma_{\zi(z_i)}^i
f_G(\zi|z_i)d\zi\right) f_{G_i}(z_i) L_t^{z_i}(W)dz_i \\ 
& =  \E \int_{-\infty}^{\infty}h_i(z_i) L_t^{z_i}(W)dz_i \\ 
& = \E \int_{0}^{t}h_i(W_s)ds, 
\end{align*} 
where in the last equality we have again used the Occupation Time 
formula. Consequently, we have
$$ \E\left[f(W_t,G)\right] = 
\frac{1}{2}\int_{0}^{t}\E\left[f_G''(W_s)\right]ds
+\frac{1}{2}\sum_{i=1}^{2d}\int_{0}^{t}\E \left[h_i(W_s)\right]ds.$$
 
If conditions 5 and 6 hold true, we have just obtained
$$|\E\left[f(W_t,G)\right]| \leq C_1t+ C_d t. $$ 

If conditions 5' and 6' hold, we apply It\^o's formula to~$f_G''(W_s)$
and get
$$\E\left[f(W_t,G)\right] = 
\frac{1}{4}\int_{0}^{t}\int_{0}^{s}\E\left[h_0''(W_u)\right]duds
+\frac{1}{4}\sum_{i=1}^{2d}\int_{0}^{t}\int_{0}^{s}\E 
\left[h_i''(W_u)\right]duds, $$ 
from which the desired result follows. 
\end{proof}

\subsection{Proof of 
Inequality~\eqref{eq:inequality-CLT-for-brownian-bridge}} 
\label{proof-key-inequality}

In this subsection we prove the 
inequality~\eqref{eq:inequality-CLT-for-brownian-bridge} which we have 
used in the proof of Lemma~\ref{lem:CLTBrownianBridge}.

Consider a family of independent random variables
$G_i\sim\mathcal{N}(0,1)$ and set
$\bar{G}^n:=n^{-1}\sum_{i=1}^{n}G_i$. 
Define the function $q(x,a,b):\R^3\to\R$ by
$$q(x,a,b) := |x-a||x-b|-|a||b|. $$ 
We aim to prove: For any $k\notin\{i,i+1\}$ one has
\begin{equation}\label{eq:bound-for-q-times-q} 
\E\left[ q(\bar{G}^n,G_{i},G_{i+1})q( \bar{G}^n,G_{k},G_{k+1}) 
\right]\leq\frac{C}{n^{2}}. 
\end{equation}

Notice that $\bar{G}^n$ can be expressed as the sum of a random 
variable independent of $(G_{i},G_{i+1},G_{k},G_{k+1})$ and a random 
variable which a.s. converges to~0: 
$$\bar{G}^n =
\frac{1}{n}\sum_{j\notin\{i,i+1,k,k+1\}}G_i + \frac{G_{i} +G_{i+1} +G_k
+G_{k+1} }{n} =: \bar{G}^{n\#}+  \frac{H_{i,k}}{n}.$$ 
In addition,
as 
$$q(x+y,a,b)\leq q(x,a,b)+ |y||2x-a-b|+|y|^2,$$
for any $y\in\R$, we have that
\begin{align*}
\E&\left[q(\bar{G}^n,G_{i},G_{i+1})q(\bar{G}^n,G_{k},G_{k+1})\right]\\
&=     \E\bigg[q(\bar{G}^{n\#}+ 
\frac{H_{i,k}}{n},G_{i},G_{i+1})q(\bar{G}^{n\#}+ 
\frac{H_{i,k}}{n},G_{k},G_{k+1})\ \bigg]  \\ 
& \leq      \E\bigg[
\left(q(\bar{G}^{n\#},G_{i},G_{i+1})+ |\frac{H_{i,k}}{n}||G_{i}+G_{i+1}
- 2\bar{G}^{n\#}|+|\frac{H_{i,k}}{n}|^2\right)\\ &\qquad\times
\left(q(\bar{G}^{n\#},G_{k},G_{k+1})+ |\frac{H_{i,k}}{n}||G_{k}+G_{k+1}
- 2\bar{G}^{n\#}|+|\frac{H_{i,k}}{n}|^2\right) \bigg]  \\ 
& =   
\E\bigg[q(\bar{G}^{n\#},G_{i},G_{i+1}) q(\bar{G}^{n\#},G_{k},G_{k+1}) 
\bigg] \\
&\quad + \E\bigg[ q(\bar{G}^{n\#},G_{i},G_{i+1})\left(
|\frac{H_{i,k}}{n}||G_{k}+G_{k+1} - 2\bar{G}^{n\#}|\right) \bigg] \\ 
& \quad+\E\bigg[ \left( |\frac{H_{i,k}}{n}||G_{i}+G_{i+1} -
2\bar{G}^{n\#}|\right) q(\bar{G}^{n\#},G_{k},G_{k+1})  \bigg]  +
\E\bigg[ \frac{H_{i,k}^2}{n^2}A \bigg], 
\end{align*} 
where $A$ is random
variable with finite second moment. An easy computation shows that
$$ \E\bigg[ \frac{H_{i,k}^2}{n^2}A \bigg] \leq \frac{C}{n^2}. $$
Since $G_{i},G_{i+1},G_{k},G_{k+1}$ are i.i.d. we have 
\begin{align*} \E\bigg[
q(\bar{G}^{n\#},G_{i},G_{i+1})&\left( |\frac{H_{i,k}}{n}||G_{k}+G_{k+1}
- 2\bar{G}^{n\#}|\right) \bigg] \\ 
&=  \E\bigg[ \left(
|\frac{H_{i,k}}{n}||G_{i}+G_{i+1} - 2\bar{G}^{n\#}|\right)
q(\bar{G}^{n\#},G_{k},G_{k+1})  \bigg].
 \end{align*} 
Hence, to prove~\eqref{eq:bound-for-q-times-q}
it is enough to prove
\begin{equation}\label{eq:bound-q-times-q-Gbar-diese}
\E\bigg[q(\bar{G}^{n\#},G_{i},G_{i+1}) q(\bar{G}^{n\#},G_{k},G_{k+1}) 
\bigg] \leq \frac{C}{n^2} 
\end{equation} 
and
\begin{equation}\label{eq:bound-for-the-rest-q-times-q} 
\E\bigg[
q(\bar{G}^{n\#},G_{i},G_{i+1})\left( |\frac{H_{i,k}}{n}||G_{k}+G_{k+1} -
2\bar{G}^{n\#}|\right) \bigg] \leq \frac{C}{n^2}. 
\end{equation} 
Notice
that $\bar{G}^{n\#}$ is independent of $(G_{i},G_{i+1},G_{k},G_{k+1})$
and has the same probability distribution as a Brownian motion $W$ at 
time~$(n-4)/n^2$. Therefore, to prove the two preceding inequalities we 
aim to use~Proposition~\ref{prop:ItoTanakaF}. This step requires easy 
calculations only. We postpone it to the Appendix.

\section{Numerical experiments}\label{sec:num-exp}

In this section we present some numerical experiments which illustrate
our main theorem~\ref{thm:TEST} and its limitations when applied to 
synthetic data. For two families of examples we examine the practical choice of the parameters~$m$ and $n$ and their effects on the Empirical Rejection Rates and Type~II errors. We also numerically compare the power of the test and the significance level in a neighborhood of the null hypothesis.

The interested reader can find an extended version 
of this section in our Supplementary Material~\cite{SupplementaryMaterial}.

\subsection{Experiments under $\mathbf{H_0}$}

In this subsection we consider the r.v.~$X = |G|^{-r}$ with $G\sim
\mathcal{N}(0,1)$ and $r>0$. Notice that $X$ has finite moments of 
order smaller than $1/r$ and that
$X\in\DA(2)$ when~$r\in(0,1/2]$.

We have performed simulations for the parameter~$r$ 
taking values in~$\{0.1, 0.2, 0.3, 0.4,0.45\}$ and for  levels 
of confidence $q=0.1$ and $q=0.05$. For each value of $r$, 
we have generated {$10^4$ independent} samples with sizes from 
$m=10^5$ up to 
$m=10^8$. In the cases $r=0.4$ and $r=0.45$ we have also considered 
$m=10^9$. To study the respective effects of $m$ and $n$, for each 
sample size $m$ we have made the {number of time steps~$n$ vary} 
from~$n=10^1$ up to $n=10^4$. 

According to Theorem~\ref{thm:TEST} one expects that the larger $m$ 
and $n$ are, the closer to $q$ the empirical rejection rate is. The 
experiments
show that increasing~$m$ improves the approximation of the 
empirical rejection rates to their expected values, whereas increasing 
$n$~does not have the same effect. One actually observes that the 
empirical rejection rates become close to~1 when $n$ becomes too large. 
Moreover, it seems that {$m$ being fixed} the optimal choice of $n$ 
depends on the value  of~$r$. 

In other words, even if our theoretical result is valid when $n$ goes to 
infinity, in practice, we should use $n$ significantly smaller than 
$m$. In addition, as $m$ increases, the convergence of the 
process $Z^m$ can be very slow. Therefore,
if $n$ is chosen too large, one is zooming in too much and sees the 
discontinuities that $Z^m$ has by construction. We continue the discussion
on that issue in Sections~\ref{sec:experiments-under-H1} and~\ref{sec:heuristic-explanation}.

In Tables 
\ref{tab:best-parameters-q=0.1} and \ref{tab:best-parameters-q=0.05} we report the experimental  values of $m$ and $n$
for which, under $\mathbf{H_0}$, the best empirical rejection rate is 
attained. Notice 
that the best value for $n$ decreases when the parameter~$r$ tends
to the critical value~$r=0.5$ (that is, when the subjacent random variable has lower and lower  finite moments). 
In particular, it does not seem to exist an optimal selection for the 
parameters $m$ and $n$ which is satisfying for any random variable in 
the domain of attraction of the Gaussian law.

\begin{table}[h]
\begin{subtable}{.47\linewidth}\centering
\begin{tabular}{|c|c|c|c|} \hline
$r$ & $m$ &$n$ & ERR \\ \hline
$0.1$ & $10^8$&$10^1,10^2$ & $0.100$ \\ \hline 
$0.2$ & $10^8$&$10^3$ & $0.100$\\ \hline 
$0.3$ & $10^7$ &$10^1$&$0.100$ \\ \hline 
$0.4$ & $10^7,10^8$&$10^1$ &$0.101$ \\ \hline 
$0.45$ & $10^8,10^9$&$10^1$ & $0.109$\\ \hline 
\end{tabular}
\caption{$q=0.1$.\label{tab:best-parameters-q=0.1}} 
\end{subtable}
\begin{subtable}{.47\linewidth}\centering
\begin{tabular}{|c|c|c|c|} \hline
$r$ & $m$&$n$ & ERR \\ \hline
$0.1$ & $10^8$& $10^4$& $0.050$\\ \hline 
$0.2$ & $10^8$& $10^4$& $0.050$\\ \hline 
$0.3$ & $10^8$& $10^1$& $0.050$\\ \hline 
$0.4$ & $10^8$& $10^1$& $0.053$\\ \hline 
$0.45$ & $10^8$& $10^1$& $0.061$\\ \hline 
\end{tabular}
\caption{$q=0.05$.\label{tab:best-parameters-q=0.05}} 
\end{subtable}
\caption{$X = |G|^{-r}$. Values of the parameters $m$ and $n$ for which the best 
empirical rejection rate (ERR) of  the $\hS_n^m$-test is attained under $\mathbf{H_0}$ for 
different values of $r$ and $q$.}
\end{table}

\begin{figure}[h]
 \centering 
\begin{subfigure}[t]{0.3\linewidth} 
 \centering 
 \includegraphics[width=1.1\linewidth, trim={0.7cm 0.7cm 0.7cm 1.2cm},clip]{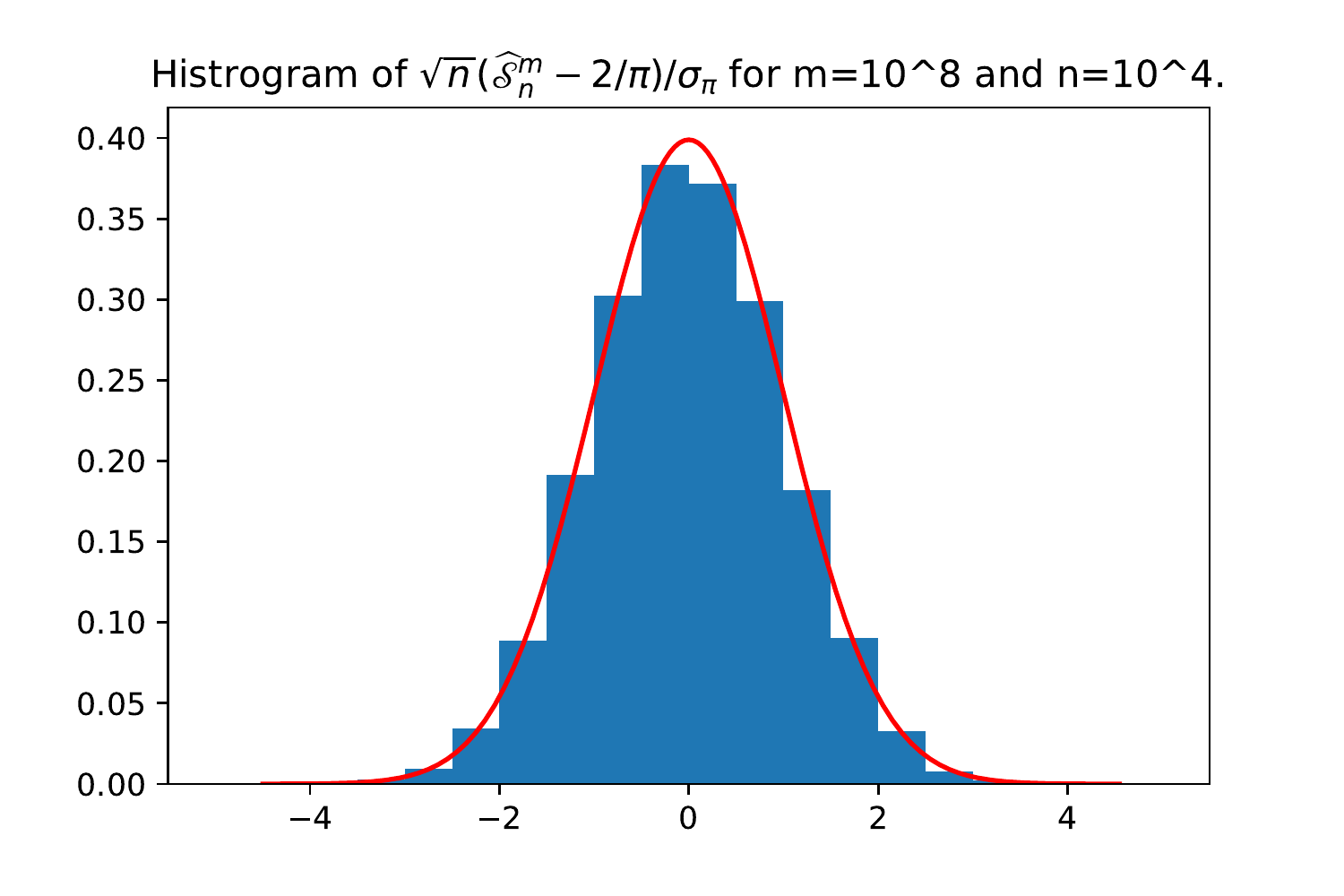}
 \end{subfigure} \quad
\begin{subfigure}[t]{0.3\linewidth} 
 \centering 
 \includegraphics[width=1.1\linewidth, trim={0.7cm 0.7cm 0.7cm 1.2cm},clip]{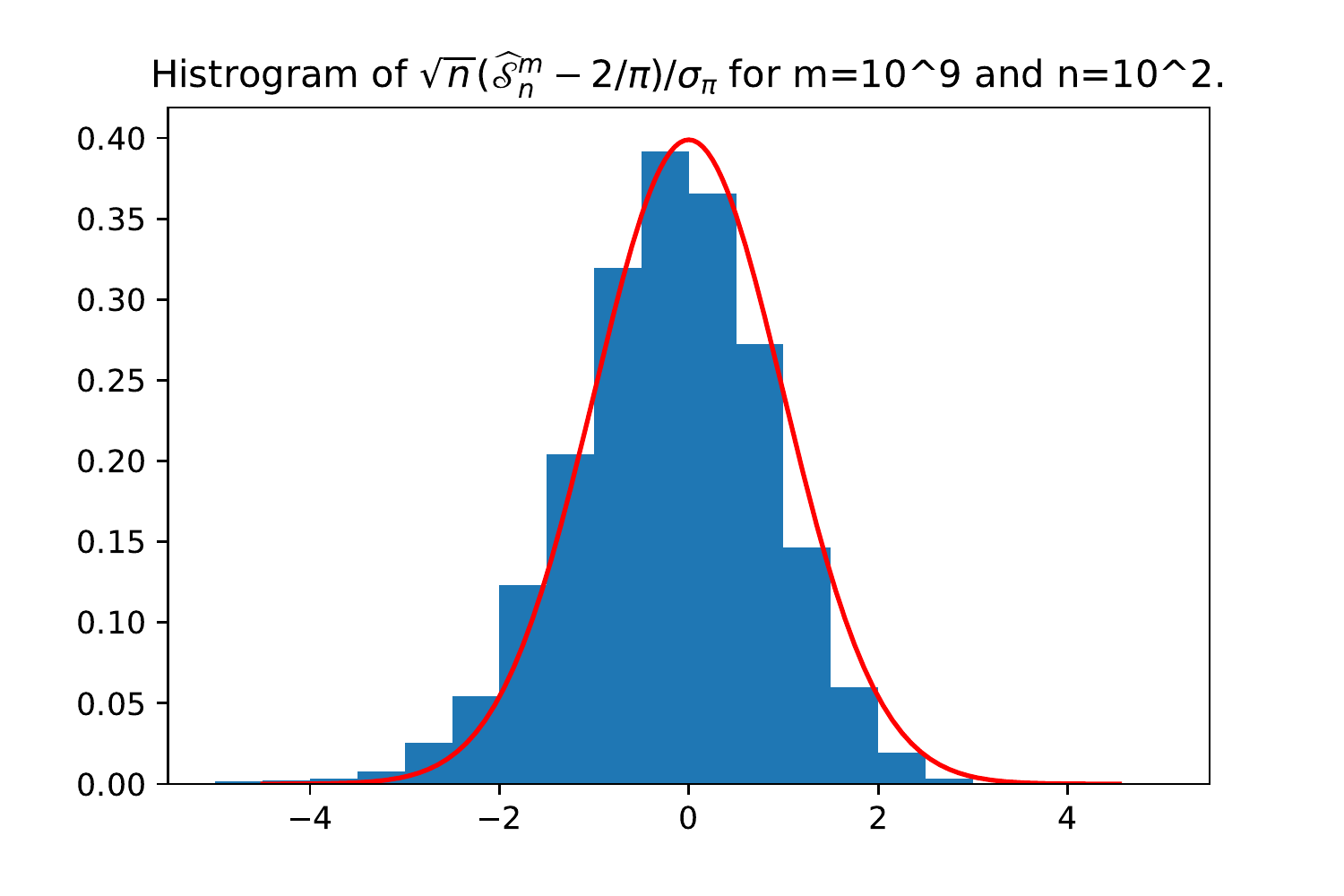}
 \end{subfigure} \quad
\begin{subfigure}[t]{0.3\linewidth} 
 \includegraphics[width=1.1\linewidth,trim={0.7cm 0.7cm 0.7cm 1.2cm},clip]{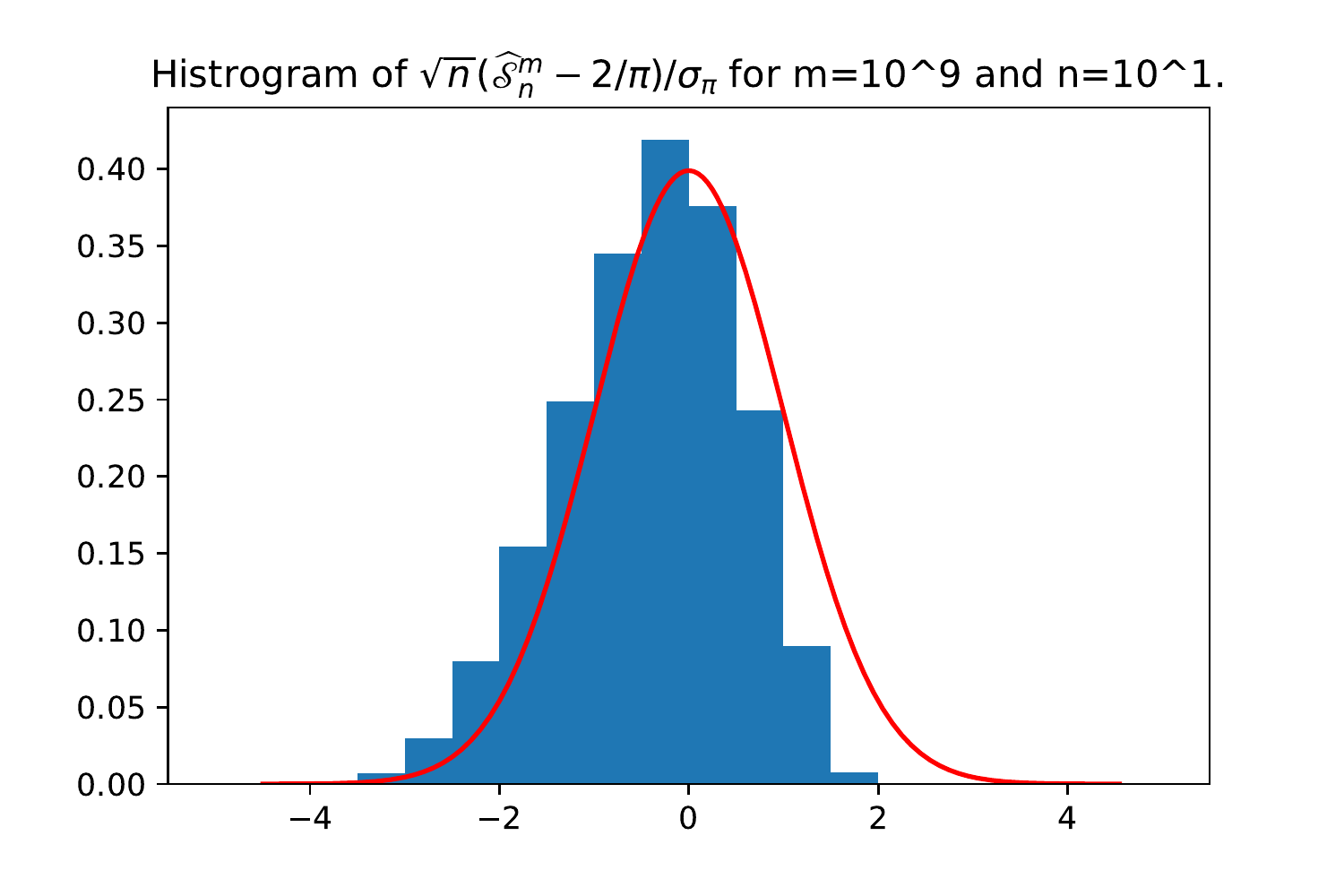}
 \end{subfigure} 
\caption{\fontsize{9}{11} $X = |G|^{-r}$.  \selectfont {The empirical distribution of} 
standardized 
${\hS_n^m}$ 
for different values of $r$, $m$ and $n$. The red line corresponds to the density of the standard normal distribution. In the left-hand side $r=0.1$, $m=10^8$ and $n=10^4$. In the center $r=0.4$,  $m=10^9$ and $n=10^2$. In the right-hand side $r=0.45$, $m=10^9$ and $n=10^2$. 
\label{figure:histogram_under_H0}} \end{figure}

In Figure \ref{figure:histogram_under_H0} appear the empirical distribution of standardized 
${\hS_n^m}$ for different values of $r$, $m$ and $n$. Notice that the adjust to the 
standard normal distribution is worse when $r$ is closer to the critical value.

\subsection{Experiments under $\mathbf{H_1}$}\label{sec:experiments-under-H1}

In this section we present empirical results which concern the Type II Error of our test, that is, the probability of failing to reject the Gaussian hypothesis $\mathbf{H_0}$ when the sample belongs to the domain of attraction $\DA(\alpha)$, with $\alpha<2$, notably when $\alpha$ is close to $2$.  Under $\mathbf{H_1}$ the empirical rejection rate (ERR) is 
expected to be close to $1$. 

We consider two cases.  First, the case where ~$X = |G|^{-r}$
with $G\sim\mathcal{N}(0,1)$  and $r>1/2$, so that $X\in\DA(1/r)$.
Second, the case where $X$  has an $\alpha$-stable distribution. 

\vspace{0.3cm}
Concerning the case where $X = |G|^{-r}$, we have performed simulations for the parameter~$r$ 
taking values in~$\{0.6, 0.75, 0.9\}$ and for  levels 
of confidence $q=0.1$ and $q=0.05$. For each value of $r$, 
we have generated {$10^4$ independent} samples with sizes from 
$m=10^5$ up to 
$m=10^8$, and to study the respective effects of $m$ and $n$, for each 
sample size $m$ we have made the {number of time steps~$n$ vary} 
from~$n=10^1$ up to $n=10^4$. To evaluate the behavior of the test  when the parameters are close to the limit case for moderate values of $n$, we also consider  
$r=0.55$ with $n\in\{200,250,500\}$. 

In Tables
\ref{tab:critical-parameters-q=0.1-h1} and \ref{tab:critical-parameters-q=0.05-h1}
we report the empirical rejection rates for different values of $r$ and 
$q$ when $n\leq 10^3$.
Notice that the rates change between $n=10^1$ and 
$n=10^2$. Tables \ref{tab:Type-II-Error-q=0.1-h1} and \ref{tab:Type-II-Error-q=0.05-h1} show that the Type II Error increases when  $r$ decreases to the limit case~$1/2$.
Tables \ref{tab:empirical-rejection-rates_01_r_0.55} and ~\ref{tab:empirical-rejection-rates_005_r_0.55}  show the satisfying behaviour of our test for moderate  values of $n$ and $m$ in the case where~$r=0.55$
for which the subjacent random variable belongs to $\DA(20/11)$.

Concerning the case where $X$ has an $\alpha$-stable distribution,  in our simulations  the parameter~$\alpha$ takes the values~$\{0.9, 1.2, 1.5, 1.8\}$. The  levels of confidence, sample sizes
and numbers of steps are as above.
 
Similarly to the preceding tables, Tables \ref{tab:ERR-q=0.1-h1-stable}, \ref{tab:ERR-q=0.05-h1-stable},
\ref{tab:ERR-q=0.05-h1-stable}, \ref{tab:ERR-q=0.05-h1-stable},
\ref{tab:Type-II-Error-q=0.1-h1-stable} and~\ref{tab:Type-II-Error-q=0.05-h1-stable} 
respectively concern empirical rejection rates and Type II Errors  for large or moderate values of~ $m$.

\begin{table}[h]
\begin{subtable}{.47\linewidth}\centering
\begin{tabular}{|c|c|c|c|} \hline
$r$ & $n=10^1$ &$n=10^2$ & $n=10^3$  \\ \hline
$0.55$&$0.157$&$0.589$&$1.000$\\ \hline 
$0.6$ & $0.194$&$0.761$ & $1.000$ \\ \hline 
$0.75$ & $0.268$&$0.950$&  $1.000$  \\ \hline 
$0.9$ & $0.340$ &$0.987$&  $1.000$  \\ \hline 
\end{tabular}
\caption{$q=0.1$.\label{tab:critical-parameters-q=0.1-h1}} 
\end{subtable}
\begin{subtable}{.47\linewidth}\centering
\begin{tabular}{|c|c|c|c|} \hline
$r$ & $n=10^1$ &$n=10^2$& $n=10^3$ \\ \hline
$0.55$&$0.102$&$0.526$&$1.000$\\ \hline 
$0.6$ & $0.123$&$0.712$&   $1.000$ \\ \hline 
$0.75$ & $0.189$&$0.934$&  $1.000$ \\ \hline 
$0.9$ & $0.251$ &$0.984$&  $1.000$ \\ \hline 
\end{tabular}
\caption{$q=0.05$.\label{tab:critical-parameters-q=0.05-h1}} 
\end{subtable}
\\ 

\begin{subtable}{.47\linewidth}\centering
\begin{tabular}{|c|c|c|c|} \hline
$r$ & $n=10^1$ &$n=10^2$ & $n=10^3$  \\ \hline
$0.55$&$0.843$&$0.411$&$0.000$\\ \hline 
$0.6$ & $0.806$&$0.239$ & $0.000$ \\ \hline 
$0.75$ & $0.732$&$0.050$&  $0.000$  \\ \hline 
$0.9$ & $0.660$ &$0.013$&  $0.000$  \\ \hline 
\end{tabular}
\caption{$q=0.1$.\label{tab:Type-II-Error-q=0.1-h1}} 
\end{subtable}
\begin{subtable}{.47\linewidth}\centering
\begin{tabular}{|c|c|c|c|} \hline
$r$ & $n=10^1$ &$n=10^2$& $n=10^3$ \\ \hline
$0.55$&$0.898$&$0.474$&$0.000$\\ \hline 
$0.6$ & $0.887$&$0.288$&   $0.000$ \\ \hline 
$0.75$ & $0.811$&$0.066$&  $0.000$ \\ \hline 
$0.9$ & $0.749$ &$0.016$&  $0.000$ \\ \hline 
\end{tabular}
\caption{$q=0.05$.\label{tab:Type-II-Error-q=0.05-h1}} 
\end{subtable}
\\

\begin{subtable}{.47\linewidth}\centering
\begin{tabular}{|c|c|c|c|} 
\cline{2-4}
\multicolumn{1}{l}{} & \multicolumn{3}{|c|}{$n$} \\ \hline
$m$&$200$&$250$&$500$\\ \hline  
$10^5$&$0.077$&$0.042$&$0.001$\\ \hline 
$10^6$&$0.126$&$0.076$&$0.004$\\ \hline 
$10^7$&$0.164$&$0.109$&$0.011$\\ \hline 
$10^8$&$0.192$&$0.130$&$0.019$\\ \hline 
\end{tabular}
\caption{$q=0.1$.\label{tab:empirical-rejection-rates_01_r_0.55}} 
\end{subtable}
\begin{subtable}{.47\linewidth}\centering
\begin{tabular}{|c|c|c|c|} 
\cline{2-4}
\multicolumn{1}{l}{} & \multicolumn{3}{|c|}{$n$} \\ \hline
$m$&$200$&$250$&$500$\\ \hline 
$10^5$&$0.105$&$0.058$&$0.002$\\ \hline 
$10^6$&$0.160$&$0.103$&$0.007$\\ \hline 
$10^7$&$0.206$&$0.139$&$0.017$\\ \hline 
$10^8$&$0.233$&$0.163$&$0.028$\\ \hline 
\end{tabular}
\caption{$q=0.05$.\label{tab:empirical-rejection-rates_005_r_0.55}} 
\end{subtable}
\caption{$X = |G|^{-r}$. \emph{First row:} Empirical rejection rates (ERR) of  the $\hS_n^m$-test
under $\mathbf{H_1}$  for $m=10^8$ and different values of $r$, $q$. \emph{Second row:} Empirical Type II Errors of  the $\hS_n^m$-test under $\mathbf{H_1}$   for $m=10^8$ and different values of $r$, $q$. \emph{Third row:} Empirical Type II Errors of the  $\hS_n^m$-test under $\mathbf{H_1}$  for $r=0.55$ and different values of $m$, $q$, $n$. \label{tab:X^G-H1} }
\end{table}

\begin{table}[h!]
\begin{subtable}{.47\linewidth}\centering
\begin{tabular}{|c|c|c|c|} \hline
$\alpha$ & $n=10^1$ &$n=10^2$ & $n=10^3$  \\ \hline
$0.9$ & $0.396$&$0.997$&$1.000$\\ \hline 
$1.2$ &$0.300$&$0.977$&$1.000$\\ \hline 
$1.5$ &$0.224$&$0.883$&$1.000$\\ \hline 
$1.8$&$0.162$&$0.611$&$1.000$\\ \hline 
\end{tabular}
\caption{$q=0.1$.\label{tab:ERR-q=0.1-h1-stable}} 
\end{subtable}
\begin{subtable}{.47\linewidth}\centering
\begin{tabular}{|c|c|c|c|} \hline
$\alpha$ &  $n=10^1$ &$n=10^2$ & $n=10^3$  \\ \hline
$0.9$ &$0.302$&$0.996$&$1.000$\\ \hline 
$1.2$ &$0.218$&$0.968$&$1.000$\\ \hline 
$1.5$ &$0.153$&$0.854$&$1.000$\\ \hline 
$1.8$&$0.102$&$0.552$&$1.000$\\ \hline 
\end{tabular}
\caption{$q=0.05$.\label{tab:ERR-q=0.05-h1-stable}} 
\end{subtable}
\\

\begin{subtable}{.47\linewidth}\centering
\begin{tabular}{|c|c|c|c|} \hline
$\alpha$ & $n=10^1$ &$n=10^2$ & $n=10^3$  \\ \hline
$0.9$ & $0.604$&$0.003$&$0.000$\\ \hline 
$1.2$ &$0.700$&$0.023$&$0.000$\\ \hline 
$1.5$ &$0.776$&$0.117$&$0.000$\\ \hline 
$1.8$&$0.838$&$0.389$&$0.000$\\ \hline 
\end{tabular}
\caption{$q=0.1$.\label{tab:Type-II-Error-q=0.1-h1-stable}} 
\end{subtable}
\begin{subtable}{.47\linewidth}\centering
\begin{tabular}{|c|c|c|c|} \hline
$\alpha$ & $n=10^1$ &$n=10^2$ & $n=10^3$ \\ \hline
$0.9$ &$0.698$&$0.004$&$0.000$\\ \hline 
$1.2$ &$0.782$&$0.032$&$0.000$\\ \hline 
$1.5$ &$0.847$&$0.146$&$0.000$\\ \hline 
$1.8$&$0.898$&$0.448$&$0.000$\\ \hline 
\end{tabular}
\caption{$q=0.05$.\label{tab:Type-II-Error-q=0.05-h1-stable}} 
\end{subtable}
\\

\begin{subtable}{.47\linewidth}\centering
\begin{tabular}{|c|c|c|c|}
\cline{2-4}
\multicolumn{1}{l}{} & \multicolumn{3}{|c|}{$n$} \\ \hline
$m$&$200$&$250$&$500$\\ \hline
$10^5$&$0.120$&$0.067$&$0.005$\\ \hline 
$10^6$&$0.149$&$0.095$&$0.010$\\ \hline 
$10^7$&$0.177$&$0.118$&$0.015$\\ \hline 
\end{tabular}
\caption{$q=0.1$.\label{tab:empirical-rejection-rates_01_alpha_1.8}}
\end{subtable}
\begin{subtable}{.47\linewidth}\centering
\begin{tabular}{|c|c|c|c|}
\cline{2-4}
\multicolumn{1}{l}{} & \multicolumn{3}{|c|}{$n$} \\ \hline
$m$&$200$&$250$&$500$\\ \hline
$10^5$&$0.154$&$0.091$&$0.007$\\ \hline 
$10^6$&$0.187$&$0.126$&$0.016$\\ \hline 
$10^7$&$0.220$&$0.150$&$0.023$\\ \hline 
\end{tabular}
\caption{$q=0.05$.\label{tab:empirical-TypeII-rates_005_alpha_1.8}}
\end{subtable}

\caption{$\alpha$-stable $X$.  \emph{First row:} Empirical rejection rates (ERR) of  the $\hS_n^m$-test
under $\mathbf{H_1}$  for $m=10^7$ and different values of $\alpha$, $q$. \emph{Second row:} Empirical Type II Errors of  the $\hS_n^m$-test
under $\mathbf{H_1}$   for $m=10^7$ and different values of $\alpha$, $q$.  \emph{Third row:} Empirical Type II Errors of the  $\hS_n^m$-test under $\mathbf{H_1}$  for $\alpha=1.8$ and different values of $m$, $q$, $n$.\label{tab:Alpha-Stable-H1}}
\end{table}

\vspace{0.3cm}
In both cases, we observe that, if one increases $m$ by letting $n$ fixed,  the quality of the test decreases. The situation seems symmetric of
what we noticed above under~$H_0$, although in this case increasing $m$ for $n$ fixed is like zooming out  the constructed trajectory
and leads the test to interpret a small jump in the limit trajectory 
as a side effect of the discretization and thus fails to reject the null hypothesis. 

\begin{figure}[h!]
 \centering 
\begin{subfigure}[t]{0.42\linewidth} 
 \centering 
  \includegraphics[width=1.1\linewidth, trim={0.7cm 0.7cm 0.7cm 0.7cm},clip]{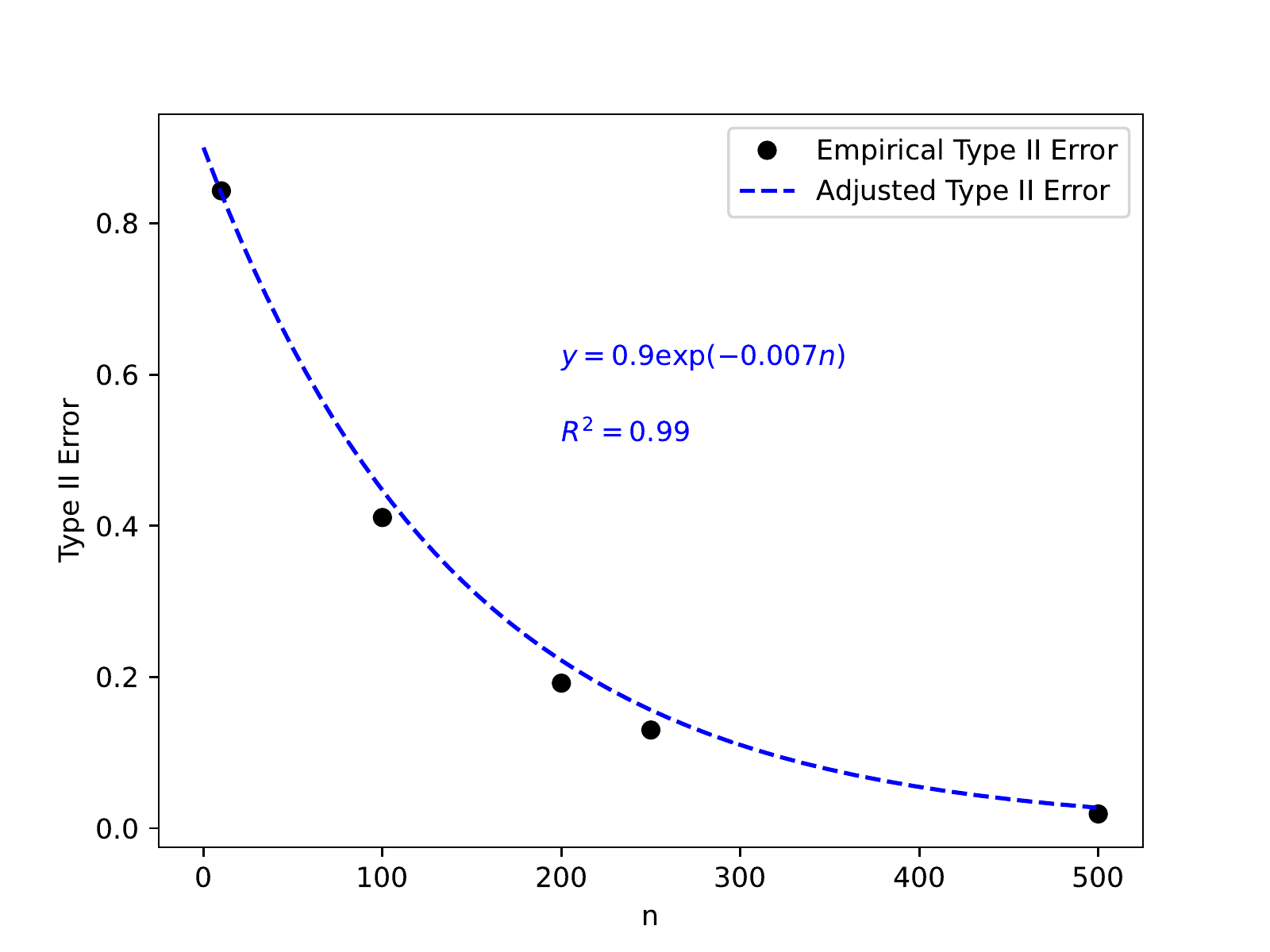}
 \end{subfigure} \quad
\begin{subfigure}[t]{0.42\linewidth} 
\centering
 \includegraphics[width=1.1\linewidth, trim={0.7cm 0.7cm 0.7cm 0.7cm},clip]{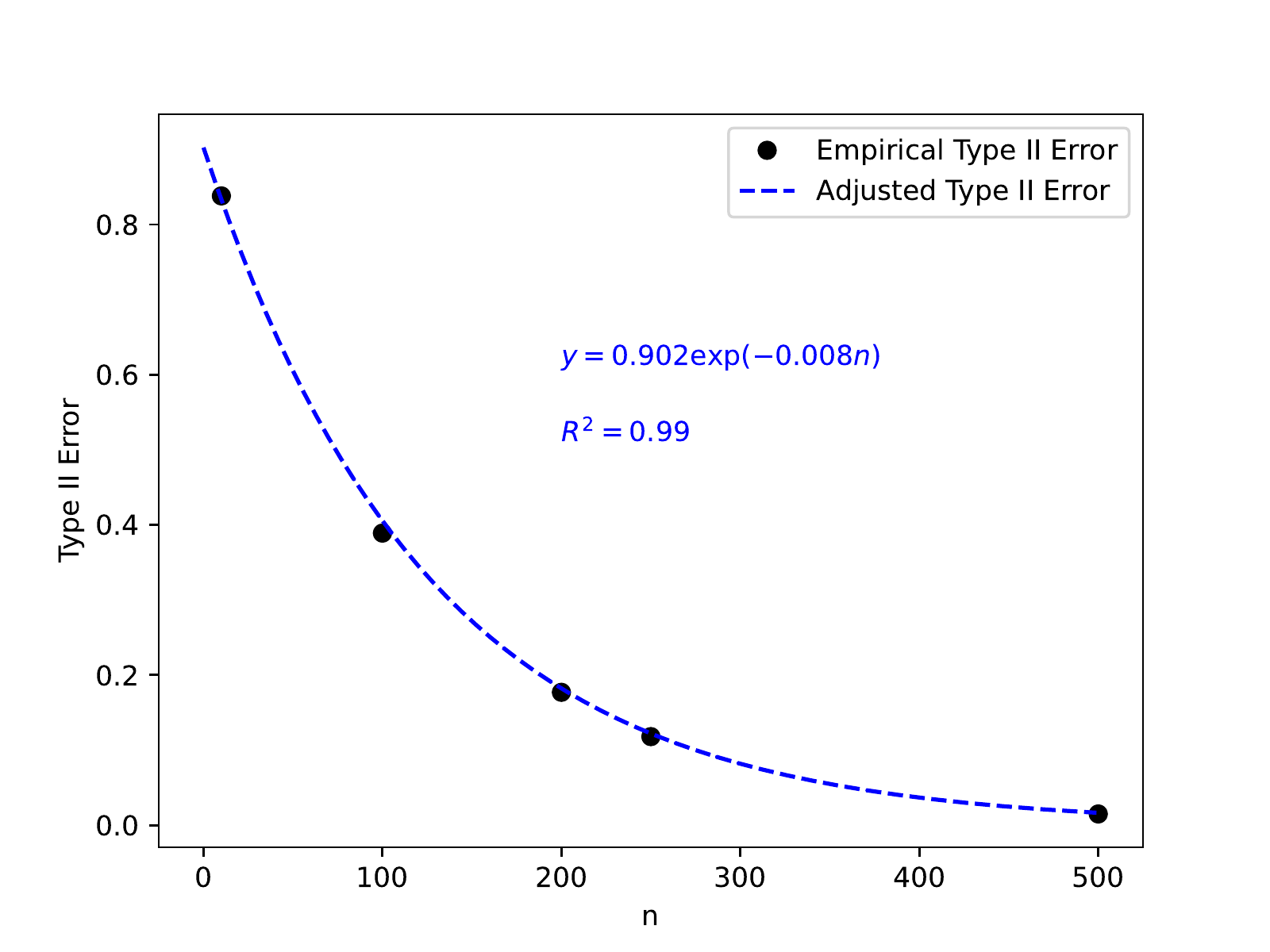}
 \end{subfigure} 
\caption{\fontsize{9}{11}   \selectfont {Empirical Type II Error as a function of $n$. \emph{Left:} $X = |G|^{-r}$ with $r=0.45$ and sample size $m=10^8$. \emph{Right:} $\alpha$-stable $X$ with $\alpha=1.8$ and  sample size $m=10^7$.  Both figures exhibit a satisfying Determination Coefficient $R^2$ larger than $0.99$.    
\label{figure:empirical-convergence-Type-II-Error}}} \end{figure}

Moreover, Tables \ref{tab:X^G-H1}  and \ref{tab:Alpha-Stable-H1} allow one to deduce the Empirical Type II Error of the Test 
$(\mathbf{ETIIE})$ for different fixed values of~$m$. We observe that it decreases, approximately, exponentially fast with $n$.
For example, when $X = |G|^{-r}$ with $r=0.45$, one observes that
$\mathbf{ETIIE}~\approx~0,9\exp\left(-0,007 n\right) $.
One thus may consider that $\mathbf{ETIIE}$ is less than $0.049$ when
$n\geq 415$. Similarly, when $X$ is $1.8$-stable, one observes that
$\mathbf{ETIIE} ~\approx~ 0,902\exp\left(-0,008 n\right)$ and one can expect that $\mathbf{ETIIE}$ is less than $0.049$ when $n\geq 364$.
See Fig.~\ref{figure:empirical-convergence-Type-II-Error}.

\vspace{0.3cm}
To summarize the tables in this section,  in the cases where $X = |G|^{-r}$ and $X$ has an
$\alpha$-stable distribution,  one can observe that for  the moderate value $n=500$ and different levels of confidence and different values of~$m$, the Type II error is much smaller than the level of confidence~$q$, even when $r$ is close to the critical value~$r=0.5$
(respectively, when $\alpha$ is close to~2).

\subsection{Some heuristic computations towards convergence rate estimates and finite-sample unbiasedness} \label{sec:heuristic-explanation}

In this subsection we discuss the important issue of the convergence rate
of our statistics when~$m$ and~$n$ tend to infinity. Accurate
convergence rate estimates would allow to prove the finite-sample
unbiasedness of our hypothesis test. As explained in
Jureckova and Kalina~\cite{jurevckova2012nonparametric},  under 
$\mathbf{H_1}$ one actually expects that the rejection probability to be at least 
the significance level of the test. See~\cite{jurevckova2012nonparametric}
for proofs of this property for several rank tests under suitable fairly general
conditions.

In our context, to accomplish this programme we have to face 
the complex algebraic expression of~$\hS_n^m$. So far, we have not 
succeeded to overcome the algebraic complexity of the necessary exact 
calculations. We thus limit ourselves to the particular case of a 
symmetric $\alpha$-stable random variable~$X$ and an analysis which is not  fully rigorous because of the use of a couple of approximations. 
However, our estimates below provide interesting insights to the respective roles of the parameters~$m$ and~$n$ and on the rate of convergence
of~$\lim_{m\to\infty}\P\left(C_{n,m} | \mathbf{H_1}\right)$ to~1 when $n$
goes to infinity. Although asymptotic, these estimates also provide
insights on the behaviour of our test for finite samples in the neighborhood of the hypothesis and  the necessary order of magnitude of the parameters~$m$ and~$n$ to ensure under $\mathbf{H_1}$ a rejection probability larger that the significance level of the test.

Recall that
\begin{equation*}
 \hS_n^m  = \frac{\sum_{i=1}^{n-1}{\left| \sum_{j=\lfloor \frac{m(i-1)}{n} 
 \rfloor +1}^{\lfloor \frac{mi}{n} \rfloor } (X_j-\overline{X}_m) 
 \right|\left| \sum_{j=\lfloor \frac{mi}{n} \rfloor +1}^{\lfloor 
 \frac{m(i+1)}{n} \rfloor } (X_j-\overline{X}_m) 
 \right|}}{\sum_{i=1}^{n}{\left|\sum_{j=\lfloor \frac{m(i-1)}{n} 
 \rfloor +1}^{\lfloor \frac{mi}{n} \rfloor } (X_j-\overline{X}_m) 
 \right|^2} }.
\end{equation*}

We are interested in cases where $\alpha\approx 2$ for which the expectation of~$|X|$ is finite . Thus the Strong Law of Large numbers applies.
Assume that $m$ is large enough to have the empirical mean $\bar{X}_m$  close to the mean of $X$, that is, to have $\bar{X}_m\approx0$~a.s., and set
\begin{equation*}
\tS_n^m  := \frac{\sum_{i=1}^{n-1}{\left| \sum_{j=\lfloor \frac{m(i-1)}{n} 
 \rfloor +1}^{\lfloor \frac{mi}{n} \rfloor } X_j 
 \right|\left| \sum_{j=\lfloor \frac{mi}{n} \rfloor +1}^{\lfloor 
 \frac{m(i+1)}{n} \rfloor } X_j 
 \right|}}{\sum_{i=1}^{n}{\left|\sum_{j=\lfloor \frac{m(i-1)}{n} 
 \rfloor +1}^{\lfloor \frac{mi}{n} \rfloor } X_j \right|^2} }.
\end{equation*}
Setting
$$ W_i := \left(\frac{n}{m}\right)^{1/\alpha} \sum_{j=\lfloor \frac{mi}{n} \rfloor +1}^{\lfloor \frac{m(i+1)}{n} \rfloor } X_j $$
we have
\begin{equation*}
\tS_n^m  = \frac{\sum_{i=1}^{n-1}|W_{i-1}||W_i|}{\sum_{i=1}^{n} W_i^2}
\end{equation*}
Notice that for any integer $i$ the random variable $W_i$ has the same distribution symmetric $\alpha$-stable as $X$.

Proceeding as in the proof of Lemma~\ref{lem:bivariation-brownian-motion-alpha-stable} we can show that the Strong Law of the Large numbers holds for the numerator. Therefore, for any~$n$ large enough,
$$
\frac{1}{n}\sum_{i=1}^{n1}|W_{i-1}||W_i| 
\approx (\E|W_1|)^2.
$$

In addition, we observe that
$$\P\left(|W_{i}|^2>t\right) = 2\P\left(W_{i}>t^{1/2}\right)
= 2\P\left(X_{i}>t^{1/2}\right)
\sim t^{\alpha/2}.$$
We thus can proceed as in the proof of Proposition \ref{prop:conv-in-P-Bivariation-Bridges} (see page 11) and apply
Theorem~4.5.2 in~\cite[Sec.5, Chap.4]{whitt2002stochastic} to
the weighted sum
$$ \frac{1}{n^{2/\alpha}} \sum_{i=1}^{n-1}|W_{i}|^2 $$
to conclude that it converges in distribution to an $\alpha/2$-stable positive random variable $\mathcal{W}$. 

Consequently, 
\begin{equation*}
 \tS_n^m  
= \frac{1}{n^{2/\alpha-1}} \frac{\frac{1}{n}\sum_{i=1}^{n-1}|W_{i}||W_{i+1}|}{\frac{1}{n^{2/\alpha}} \sum_{i=1}^{n-1}|W_{i}|^2}\\
\approx \frac{1}{n^{2/\alpha-1}} \frac{(\E|W_1|)^2}{\mathcal{W}},
\end{equation*}
where $\mathcal{W}$ is $\alpha/2$-stable. It comes:
\begin{align*}
\P\left( \hS_n^m <  \frac{2}{\pi} - z_{1-q/2}
\sqrt{\frac{\sigma_\pi^2}{n}} \right)  
&\approx 
\P\left( \tS_n^m <  \frac{2}{\pi} - z_{1-q/2}
\sqrt{\frac{\sigma_\pi^2}{n}} \right)  \\
&\approx \P\left(  \frac{1}{n^{2/\alpha-1}} \frac{(\E|W_1|)^2}{\mathcal{W}} <  \frac{2}{\pi} - z_{1-q/2}
\sqrt{\frac{\sigma_\pi^2}{n}} \right) \\
&\approx \P\left( \mathcal{W} >  \frac{1}{n^{2/\alpha-1}} \frac{(\E|W_1|)^2}{\frac{2}{\pi} - z_{1-q/2}
\sqrt{\frac{\sigma_\pi^2}{n}}} \right)\\
&\approx 1 - \P\left( \mathcal{W} \leq  \frac{1}{n^{2/\alpha-1}} \frac{(\E|W_1|)^2}{\frac{2}{\pi} - z_{1-q/2}
\sqrt{\frac{\sigma_\pi^2}{n}}} \right).
\end{align*}
Similarly,
\begin{align*}
\P\left( \hS_n^m >  \frac{2}{\pi} + z_{1-q/2}
\sqrt{\frac{\sigma_\pi^2}{n}} \right) &\approx \P\left(  \frac{1}{n^{2/\alpha-1}} \frac{(\E|W_1|)^2}{\mathcal{W}} 
>  \frac{2}{\pi} + z_{1-q/2}
\sqrt{\frac{\sigma_\pi^2}{n}} \right) \\
&\approx \P\left( \mathcal{W} <  \frac{1}{n^{2/\alpha-1}} \frac{(\E|W_1|)^2}{\frac{2}{\pi} + z_{1-q/2}
\sqrt{\frac{\sigma_\pi^2}{n}}} \right).
\end{align*}
We now use a bound for the density of $\mathcal{W}$. Let $\xi$ be
a $\overline{\alpha}$-stable random variable. Its characteristic function is of the form
$$ \exp(i~a~\lambda - c^{\overline{\alpha}}~|\lambda|^{\overline{\alpha}}
(1-i~\beta~\text{sign}(\lambda)~w(\lambda,\overline{\alpha}))) $$
for some parameters $0<\overline{\alpha}\leq 2$, $c>0$
and $-1\leq\beta\leq1$, where
\begin{equation*}
w(\lambda,\overline{\alpha}) =
\begin{cases}
\tan(\frac{\pi\overline{\alpha}}{2})~~~\text{if}~~~\overline{\alpha}\neq1, \\
-\frac{2}{\pi}~\log(|\lambda|)~~~\text{if}~~~\overline{\alpha}=1.
\end{cases}
\end{equation*}
In our case, $\mathcal{W}$ is positive and thus $\beta=1$
(see e.g. Nolan~\cite[Lem.1.1, Chap.1, Sec.1.4] {nolan2020univariate}).

We now use that the density of any $\alpha$-stable random variable with skewness parameter $\beta=1$ and scale parameter  $c$ is bounded by $\Gamma(1+1/\bar{\alpha})/c\pi$  (see Lemma 3.11 in \cite{nolan2020univariate}), so
$$
 \P\left( \mathcal{W} <  \frac{1}{n^{2/\alpha-1}} \frac{(\E|W_1|)^2}{\frac{2}{\pi} + z_{1-q/2}
\sqrt{\frac{\sigma_\pi^2}{n}}} \right) \leq  \frac{\Gamma(1+2/\alpha)}{c\pi} \frac{1}{n^{2/\alpha-1}} \frac{(\E|W_1|)^2}{\frac{2}{\pi} + z_{1-q/2}
\sqrt{\frac{\sigma_\pi^2}{n}}},
$$
$$
\P\left( \mathcal{W} \leq  \frac{1}{n^{2/\alpha-1}} \frac{(\E|W_1|)^2}{\frac{2}{\pi} - z_{1-q/2}
\sqrt{\frac{\sigma_\pi^2}{n}}} \right)\leq  \frac{\Gamma(1+2/\alpha)}{c\pi}\frac{1}{n^{2/\alpha-1}} \frac{(\E|W_1|)^2}{\frac{2}{\pi} - z_{1-q/2}
\sqrt{\frac{\sigma_\pi^2}{n}}}.
$$
Finally,
\begin{equation} \label{ineq:estimate-P-wrt-n}
\begin{split}
\P(C_{n,m}|\mathbf{H}_1) &=  \P\left( \hS_n^m <  \frac{2}{\pi} - z_{1-q/2}
\sqrt{\frac{\sigma_\pi^2}{n}} \right) + \P\left( \hS_n^m >  \frac{2}{\pi} + z_{1-q/2}
\sqrt{\frac{\sigma_\pi^2}{n}} \right) \\
& \geq 1 - \frac{C(\alpha,q)}{ n^{\frac{2}{\alpha}-1}} ,
\end{split}
\end{equation}
where the constant $C(\alpha,q)$ remains uniformly bounded for $\alpha\in(1,2)$.
Therefore, up to the couple of approximations made in the above calculation and the fact that $\alpha\in(1,2)$, $\P(C_{n,m}|\mathbf{H}_1)$ increases to~1  at a rate close to
$n^{2/\alpha-1}$.
Notice that the parameter $m$ is needed large enough only to ensure
that $\bar{X}_m\approx0$~a.s. The following useful convergence rate estimate w.r.t.~$m$ 
can be found in Basu et al.~\cite{Basu-al}:

\begin{thm}
Let $(X_j)$ be a sequence of i.i.d. randonm variables with common probability density~$v_1$. Suppose that they are centered and belong to the domain
of attraction of a stable law of index~$1<\alpha<2$ whose probability density is denoted by $v_\alpha$. Suppose that their
characteristic function belongs to $L^r(\R)$ with~$r\geq 1$. Finally,
suppose that $\int_{\R} x^2~|v_1(x) - v_\alpha(x)|~dx < \infty$. Then, for some positive number~$A$, for any~$m$ large enough, the probability density $v_m$ of the normalized sum
$Z_m := \frac{A}{m^{1/\alpha}}~\sum_{i=1}^m X_i$ satisfies:
$$ \sup_x (1+|x|^\alpha)~|v_m(x)-v_\alpha(x)| = 
\mathcal{O}(\frac{1}{m^{\frac{2}{\alpha}-1 } }) . $$ 
\end{thm}
In the case of our simulations where the common law of the $X_j$'s  is stable
we deduce that
\begin{equation*}
\begin{split}
\P(|\overline{X}_m| > \epsilon) &= \P(|Z_m| > A\epsilon~m^{1-\frac{1}{\alpha}} ) \\
&\leq \int_{|x|>A\epsilon m^{1-\frac{1}{\alpha}}} v_\alpha(x)~dx
+ \mathcal{O}(\frac{1}{m^{\frac{2}{\alpha}-1}}),
\end{split}
\end{equation*}
We now use tail estimates from Gairing and Imkeller~\cite[Sec.4.1]{Imkeller}. 
Let $\xi$ be a centered $\overline{\alpha}$-stable random variable with  characteristic function as above.  Then,
For any $x>0$ it holds that
$$ \P(\xi>x) = \frac{1}{\pi}~\frac{\overline{\alpha}(\overline{\alpha}+1)}{\overline{\alpha}}
~\sin((\overline{\alpha}+\delta)\tfrac{\pi}{2})
~\frac{c^{\overline{\alpha}}}{\cos(\delta\tfrac{\pi}{2})}
~x^{-\overline{\alpha}}  + \mathcal{O}(x^{-2\overline{\alpha}}) $$
with $\delta := \frac{2}{\pi}~\arctan(\beta \tan(\overline{\alpha}\tfrac{\pi}{2}))$. 
Similarly,
$$ \P(\xi < -x) = \frac{1}{\pi}~\frac{{\overline{\alpha}}(\overline{\alpha}+1)}
{\overline{\alpha}}~\sin((\overline{\alpha}-\delta)\tfrac{\pi}{2})
~\frac{c^{\overline{\alpha}}}{\cos(\delta\tfrac{\pi}{2})}
~x^{-\overline{\alpha}} + \mathcal{O}(x^{-2\overline{\alpha}}). $$
Here, $\overline{\alpha}=\alpha$ and we deduce:
\begin{equation} \label{ineq:cv-rate-m}
\P(|\overline{X}_m| > \epsilon)
\leq \frac{C}{\epsilon^{\alpha}~m^{\alpha-1}} 
+ \frac{C}{m^{\frac{2}{\alpha}-1}}.
\end{equation}

In the preceding calculations,  the approximation of $\hS_n^m$ by
the more tractable quantity~$\tS_n^m$ is based on the hypothesis 
that the common law of the $X_j$'s  is stable and symmetric.
In contrast, the result in \cite{Basu-al} do not suppose 
that the $X_i$'s have a stable law: They actually apply to the case of random
variables in the domain of attraction of a stable law. Anyway, our estimates
are certainly sub-optimal.
For example, for $\alpha=1.8$ the heuristic convergence rate is $n^{-0.1}$ whereas  the results in~Table \ref{tab:Alpha-Stable-H1}  tend to show
an exponential decay.

However, whilst not perfect, these estimates give some insight on the reason for which~$m$ needs
to be chosen much larger than~$n$, especially when $\alpha$ is close to $2$.
Actually, $\overline{X}_m$ can be seen as a random perturbation term
in $\hS_n^m$ whose
expectation needs to be much smaller than
the right-hand side of~\eqref{ineq:estimate-P-wrt-n}.

Finally, notice that, up to the approximation of $\bar{X}_m$ made in the above calculation, the inequalities~\eqref{ineq:estimate-P-wrt-n}
and~\eqref{ineq:cv-rate-m} show that, for any small $\epsilon$ and  any~$m$ and~$n$ large enough
we have $\P(C_{n,m}|\mathbf{H}_1)>1-\epsilon$. which suggests that
our test is finite-sample unbiased.

\subsection{Numerical experiments under weak dependence} \label{sec:num-exp-weak-dep}

The objective of this section is to show we can extend
our test to some cases of weakly dependent data. Certainly, extra 
assumptions are necessary to obtain functional limit theorems by using 
the $J_1$-topology as in Theorem~\ref{thm:donsker}: 
For example, Avram and Taqqu~\cite{10.1214/aop/1176989938} show that 
the $J_1$-weak convergence of normalized sums cannot hold for the 
classes of dependent 
random variables they consider.

Let us exhibit some sufficient conditions to replace 
the independency condition in  Theorem~\ref{thm:TEST}.  Assume that the 
sample is stationary, $r$-dependent
and satisfy the two following conditions:

\begin{equation}\label{eq:cond-mixing-Shao}
\sup\left\{ |\Cor(f,g)|:\real f\in L^2(\sigma(X_1)), \real g\in L^2(\sigma(X_{2},X_{3},\ldots)) \right\} <1
\end{equation}
and
\begin{equation}\label{eq:cond-m-dependent}
\left(\forall \epsilon>0,\, \forall j=2,\ldots,k \right)\;\;
\lim_{m\to\infty}\P\left(|X_{j}|>\epsilon c_m/|X_1|>\epsilon c_m \right)=0.
\end{equation}

Let us start with adapting the proof of Proposition \ref{prop:conv-in-P-estimator} to this setting. First, if $X$ is in the domain of attraction of the normal law, in view of \eqref{eq:cond-mixing-Shao} Corollary 1.1 in Shao~\cite{Qiman:1993aa}  implies that $(S_{\lfloor mt 
\rfloor}-\mu_n)/c_m$ converges to a scaled Brownian motion, from which  $Z^m$ converges to a Brownian bridge. Second, 
if $X$ belongs to the domain of attraction of a stable law, under  \eqref{eq:cond-m-dependent}  Corollary~1.4 in Tyran-Kamińska \cite{TYRANKAMINSKA20101629} states 
that $(S_{\lfloor mt \rfloor}-\mu_m)/c_m$ converges in distribution to a 
$\alpha$-stable random process, from which $Z^m$ converges in 
distribution to a $\alpha$-stable bridge. The rest of the arguments in the proof
of Proposition  \ref{prop:conv-in-P-estimator} do not rely on the independency hypothesis. 

Let us now turn to Proposition \ref{prop:central-limit-theorem-statistic}. Assume that $X$ is in the domain of attraction of the normal law. In view of \eqref{eq:cond-mixing-Shao}  $Z^m$ converges to a Brownian bridge.
Therefore, under $\mathbf{H_0}$,
$$ \left|  
\E\left[\psi\left(\frac{\sqrt{n}}{\sigma_\pi}
({\hS_n^m}-\frac{2}{\pi})\right)\right]
- \E\left[\psi\left(\frac{\sqrt{n}}{\sigma_\pi}(\hS_n(Z)
-\frac{2}{\pi})\right)\right]\right|, $$
converges to zero,  for every fixed $n$,  as $m$ goes to infinity. All the other arguments in the proof of Proposition 
\ref{prop:central-limit-theorem-statistic} remain valid in this weak dependence setting. 

We thus have shown one possible extension for Theorem~\ref{thm:TEST} to samples of non-independent random variables.

To illustrate the behavior of the test for non-independent random variables, 
we consider a sequence 
$Y_2,\ldots,Y_m$ of 
$1$-dependent random variables constructed as follows. First, we 
generate a sample $X_1,\ldots,X_m$ of independent copies
of ~$|G|^{-r}$ where $G\sim 
\mathcal{N}(0,1)$, and then
\begin{equation}\label{eq:def_weak_dependent}
Y_k = \frac{X_{k-1}}{X_{k-1}+1}X_k.
\end{equation}

Second, we apply our methodology to $Y_2,\ldots,Y_m$ for $r=0.3$ and 
$r=0.75$. Just as above we consider $10^4$ scenarios. 

For $r=0.3$, $\mathbf{H_0}$ holds and the experiment shows 
that as $m$ increases the empirical rejections rates tend 
to the theoretical ones, for example for $q=0.1$, $m=10^8$ and $n=10^3$, 
the empirical rejection rate is $0.11$. As in the i.i.d. case we also observe that
the performance of the methodology decreases when~$n$ becomes too 
large.

For $r=0.75$, that is, under 
$\mathbf{H_1}$, smaller values of $n$ 
leads to bigger `type II' error. Nevertheless, for $n=10^2$ we 
observe empirical rejection rates close to~1.

\begin{figure}[h] 
 \centering 
 \includegraphics[width=0.5\linewidth, trim={0.7cm 0.7cm 0.7cm 1.2cm},clip]{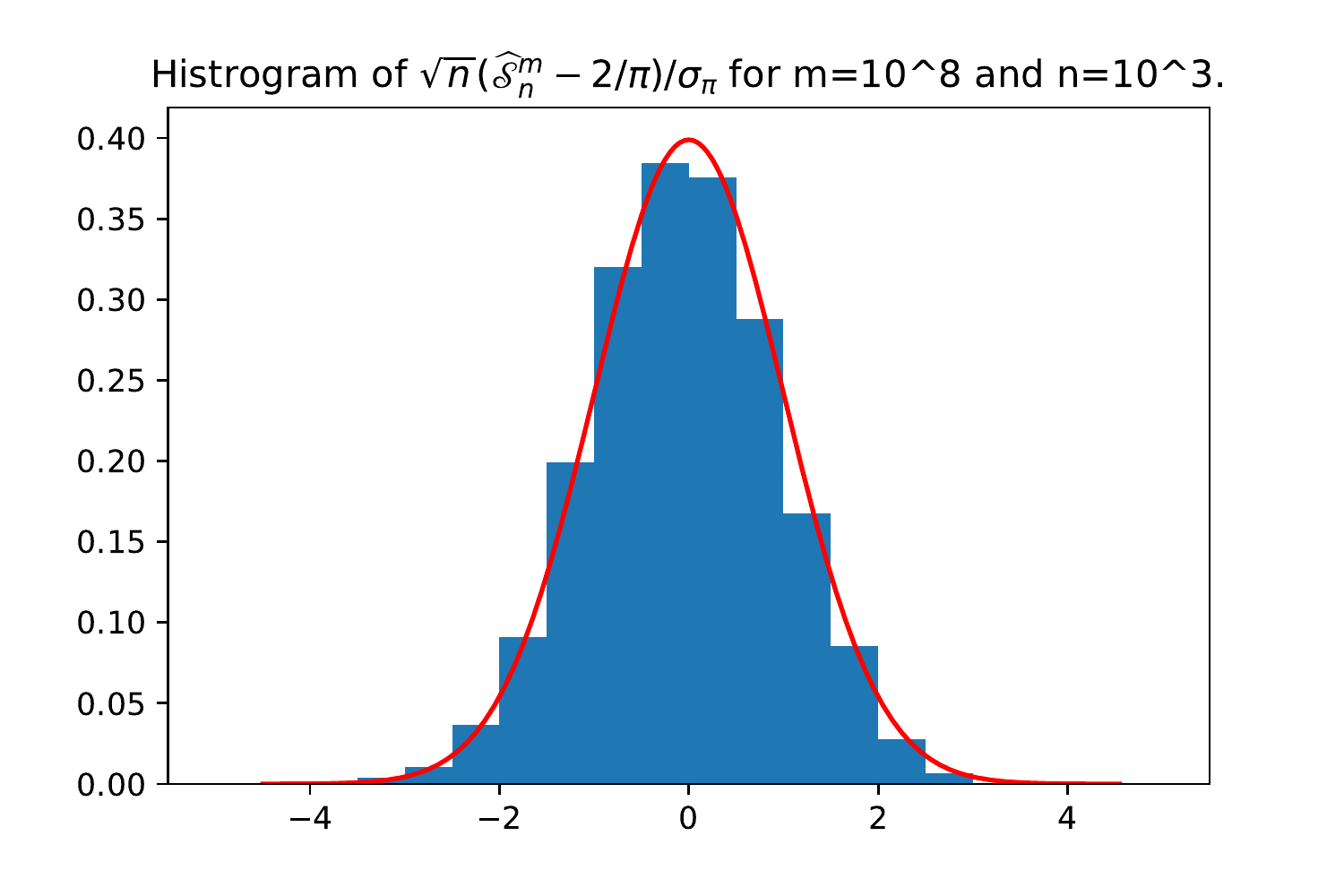}
\caption{\fontsize{9}{11}\selectfont {Empirical distribution of} 
standardized 
${\hS_n^m}$ 
for the variables $Y_k$ defined in \eqref{eq:def_weak_dependent}.  The red line corresponds to the density of the standard normal distribution.
\label{figure:histogram_under_H0_weak}} \end{figure}

\section{Conclusion}

Inspired by Barndorff-Nielsen and Shephard's methodology 
to test the existence of jumps in the trajectories of certain 
semimartingales, we have built an statistical test to determine if a 
given random variable $X$ belongs to the domain of attraction of the 
Gaussian law or of another stable law.
This statistical test allows us to give a partial answer to the 
ill-posed (but important in practice) problem of testing
the finiteness of moments of an observed random variable.

So far, our methodology is based on asymptotic results. For
the purpose of practical applications it would be interesting to obtain 
a non-asymptotic version of our test or, at least, to obtain 
theoretical results on suitable dependences 
between the parameters $m$ and $n$ which guarantee that our hypothesis 
test is reliable. {To obtain accurate estimates w.r.t. $m$ and $n$
for $\P\left(C_{n,m}|H_0\right)$ and $\P\left(C_{n,m} | H_1\right)=1$
(which supposes to obtain non asymptotic estimates precising 
the CLT for~$\hS_n^m$) the calculations appear to be heavy, lengthy and 
very technical. We hope to be able to obtain satisfying 
results in the next future.}\\

\noindent\textbf{Acknowledgment: } The two authors thank the referees 
for their useful suggestions to improve the overall organisation of the 
manuscript and to clarify some mathematical issues.

The first author is also grateful to Karine 
Bertin and Joaquin Fontbona for some helpfull discussions and comments 
on a previous version of this work.

\appendix\label{sec:appendix} 
\appendixpage %

\section{An elementary lemma on Gaussian distributions}
In the proof below of the 
inequalities~\eqref{eq:bound-q-times-q-Gbar-diese} 
and \eqref{eq:bound-for-the-rest-q-times-q}
we use the following elementary property of standard Gaussian laws.
Recall that at the beginning of Section~\ref{proof-key-inequality} the 
function~$q$ was defined as
$$ q(x,a,b)=|x-a||x-b|-|a||b|. $$

Let $\phi$, $\Phi$ respectively denote the
density and the cumulative distribution function of a standard Gaussian
random variable. In the computations below we will use the following elementary equalities
that hold for any $x\in\R$ :
\begin{equation}\label{eq:sign-G-x}
\int_{-\infty}^{\infty}\sgn(x-z_1)~\phi(z_1)~dz_1 =  2\Phi(x) -1,
\end{equation} 
\begin{equation}\label{eq:abs-G-x}
\int_{-\infty}^{\infty}|x-z_2|~\phi(z_2)~dz_2 
=   2x\Phi(x) - x  + 2\phi(x),
\end{equation} 
and 
\begin{equation}\label{eq:q-G-x}
\int_{-\infty}^{\infty}\int_{-\infty}^{\infty}
q(x,z_3,z_4)~\phi(z_3)~\phi(z_4)~dz_3~dz_4 =  \left(2x\Phi(x) - x  
+ 2\phi(x)\right)^2-\frac{2}{\pi}. 
\end{equation}

\subsection{Proof of Inequality~\eqref{eq:bound-q-times-q-Gbar-diese}} 
Let $f(x,z_1,z_2,z_3,z_4)$ be defined as
$$ f(x,z_1,z_2,z_3,z_4) := q(x,z_1,z_2)q(x,z_3,z_4). $$ 
Then
$$D_{-}f(x,z_1,z_2,z_3,z_4)=D_{-}q(x,z_1,z_2)q(x,z_3,z_4)
+q(x,z_1,z_2)D_{-}q(x,z_3,z_4),$$
and 
\begin{align*} 
D^2f(dx,z_1,z_2,z_3,z_4) &=
f''(x,z_1,z_2,z_3,z_4)dx\\ &\quad+ 
2|z_1-z_2|q(z_1,z_3,z_4)\delta_{z_1}(dx)
+2|z_2-z_1|q(z_2,z_3,z_4)\delta_{z_2}(dx)\\
&\quad+2|z_3-z_4|q(z_3,z_1,z_2)\delta_{z_3}(dx)+2|z_4-z_3|
q(z_4,z_1,z_2)\delta_{z_4}(dx),
\end{align*} 
where 
\begin{align*} 
f''(x,z_1,z_2,z_3,z_4) &=  
2\left[\sgn(x-z_1)|x-z_2|+|x-z_1|\sgn(x-z_2)\right]\\ 
&\quad\times\left[\sgn(x-z_3)|x-z_4|+|x-z_3|\sgn(x-z_4) \right] \\ 
&\quad+2\sgn(x-z_1)\sgn(x-z_2) q(x,z_3,z_4)\\ &\quad 
+2\sgn(x-z_3)\sgn(x-z_4) q(x,z_1,z_2). 
\end{align*} 
It is then easy to check that the conditions 1-4 4 of 
Proposition~\ref{prop:ItoTanakaF} hold true.

It remains to check that conditions 5' and 6' of Proposition
\ref{prop:ItoTanakaF} also hold true. We use the same notation as in 
the aforementioned proposition. Observe that
\begin{align*} 
h_0(x)&=
\E\left[f''(x,G_{i},G_{i+1},G_{k},G_{k+1})\right] \\ &=
\int_{-\infty}^{\infty}\int_{-\infty}^{\infty} \int_{-\infty}^{\infty}
\int_{-\infty}^{\infty}f''(x,z_1,z_2,z_3,z_4)
\phi(z_1)\phi(z_2)\phi(z_3)\phi(z_4)dz_1dz_2dz_3dz_4 \\ 
& = 8\left(\int_{-\infty}^{\infty}\sgn(x-z_1)\phi(z_1)dz_1\right)^2
\left(\int_{-\infty}^{\infty}|x-z_2|\phi(z_2)dz_2\right)^2\\
&\quad+4\left(\int_{-\infty}^{\infty}\sgn(x-z_1)\phi(z_1)dz_1\right)^2
\int_{-\infty}^{\infty}\int_{-\infty}^{\infty}
q(x,z_3,z_4)\phi(z_3)\phi(z_4)dz_3dz_4. 
\end{align*} 
In view of equalities  \eqref{eq:sign-G-x}, \eqref{eq:abs-G-x} and \eqref{eq:q-G-x}
a straightforward computation leads to
\begin{equation*}\label{eq:} 
\begin{aligned} 
   h_0(x) &= 12\left( 2\Phi(x) -1\right)^2\left( 2x\Phi(x) - x  + 
   2\phi(x)\right)^2-\frac{8}{\pi}\left( 2\Phi(x) -1\right)^2. 
\end{aligned} 
\end{equation*} 

It is clear that $h_0(0)=0$ and also 
that $h$ is smooth enough to satisfy Condition 5' in 
Proposition~\ref{prop:ItoTanakaF}.

On the other hand, the weights of the singular part of
$D^2f(dx,z_1,z_2,z_3,z_4)$ are given by 
\begin{align*} 
\gamma^1_z &= 
2|z_1-z_2|q(z_1,z_3,z_4),& \gamma^2_z &= 2|z_2-z_1|q(z_2,z_3,z_4),\\
\gamma^3_z &=  2|z_3-z_4|q(z_3,z_1,z_2),& \gamma^4_z &=
2|z_4-z_3|q(z_4,z_1,z_2). \end{align*} Since $G_i,G_{i+1},G_{k},G_{k+1}$
are i.i.d. it follows \begin{align*} h_1(x) &= 
\left(\int_{-\infty}^{\infty}\int_{-\infty}^{\infty}
\int_{-\infty}^{\infty}|x-z_2|q(x,z_3,z_4)\phi(z_2)
\phi(z_3)\phi(z_4)dz_2dz_3dz_4\right)
\phi(x)\\ 
&=\left(\int_{-\infty}^{\infty}|x-z_2|\phi(z_2)dz_2\right)
\left(\int_{-\infty}^{\infty}\int_{-\infty}^{\infty}
 q(x,z_3,z_4)\phi(z_3)\phi(z_4)dz_3dz_4\right)\phi(x)      \\ 
& =   \left(2x\Phi(x) - x  + 2\phi(x)\right)\left(\left(2x\Phi(x) - x  
+ 2\phi(x)\right)^2-\frac{2}{\pi}\right)\phi(x) . 
\end{align*} 
Notice that $h_1$ is $\mathcal{C}^\infty$, and its derivatives have 
polynomial growth, therefore $h_1$ satisfies Condition 6' in 
Proposition~\ref{prop:ItoTanakaF}. 

Since $G_i,G_{i+1},G_{k},G_{k+1}$ 
are identically distributed, one has that 
$h_2(x)=h_3(x)=h_4(x)=h_1(x)$, 
and therefore $h_2,h_3$ and $h_4$ also satisfy  Condition 6' in 
Proposition \ref{prop:ItoTanakaF}. 

\section{Proof of Inequality~\eqref{eq:bound-for-the-rest-q-times-q}} 

Notice that 
\begin{align*} 
\E&\bigg[ q(\bar{G}^{n\#},G_{i},G_{i+1})\left(
|\frac{H_{i,k}}{n}||G_{k}+G_{k+1} - 2\bar{G}^{n\#}|\right) \bigg]    \\
& \leq    \E\bigg[ q(\bar{G}^{n\#},G_{i},G_{i+1})\left(
|\frac{G_{i}+G_{i+1}}{n}|+
|\frac{G_{k}+G_{k+1}}{n}|\right)\left(|G_{k}+G_{k+1}| +
2|\bar{G}^{n\#}|\right) \bigg]     \\ 
& =     \E\bigg[
q(\bar{G}^{n\#},G_{i},G_{i+1}) |\frac{G_{i}+G_{i+1}}{n}||G_{k}+G_{k+1}|
\bigg]   +   2\E\bigg[ q(\bar{G}^{n\#},G_{i},G_{i+1})
|\frac{G_{i}+G_{i+1}}{n}||\bar{G}^{n\#}| \bigg]   \\ 
&\quad +  \E\bigg[
q(\bar{G}^{n\#},G_{i},G_{i+1}) \frac{|G_{k}+G_{k+1}|^2}{n} \bigg]   +  2
\E\bigg[ q(\bar{G}^{n\#},G_{i},G_{i+1})
|\frac{G_{k}+G_{k+1}}{n}||\bar{G}^{n\#}| \bigg]\\ 
&=: E_1+E_2+E_3+E_4.
\end{align*} 
We are now going to bound from above each $E_i$. 

Again notice
that $\bar{G}^{n\#}$  is equal in law to a Brownian motion $W$ at time
$(n-4)/n^2$, from which
\begin{align*} 
E_2 &=  \frac{2}{n}\E\bigg[ 
|G_{i}+G_{i+1}|q(W_{\frac{n-4}{n^2}},G_{i},G_{i+1})|W_{\frac{n-4}{n^2}}|
\bigg], \end{align*} 
We aim to apply the proposition~\ref{prop:ItoTanakaF} to the function
$$ f(x,z_1,z_2) := |z_1+z_2|~q(x,z_1,z_2)~|x|. $$
In this case, $z=(z_1,z_2)$. We have: 
$$D_{-}f_{z}(x)
=|z_1+z_2|D_{-}q(x,z_1,z_2)|x|+|z_1+z_2|q(x,z_1,z_2)\sgn(x),$$
and 
\begin{align*} 
D^2f_{z}(dx) &= f_{z}''(x)dx+ 
2|z_1+z_2||z_1-z_2|\delta_{z_1}(dx) 
+ 2|z_1+z_2||z_2-z_1|\delta_{z_2}(dx), 
\end{align*} 
where \begin{align*}
f_{z}''(x) &=
2|z_1+z_2|\left(\sgn(x-z_1)|x-z_2|+|x-z_1|\sgn(x-z_2)\right)\sgn(x)\\
&\quad+2|z_1+z_2|\sgn(x-z_1)\sgn(x-z_2). 
\end{align*} 

We now show that for $i=0,1,2$, the functions 
$t\to \E\left[h_i(W_t)\right]$ are bounded for $t$ in a neighborhood of 
~$0$.

Indeed, \begin{align*} |h_0(x)| &= |\E\left[f_{G}''(x)\right]|     \\ 
&\leq \int_{\R^2} \left(2|z_1+z_2|\left(|x-z_2|+|x-z_1|\right)
+2|z_1+z_2|\right)\phi(z_1)\phi(z_2)dz_1dz_2 \\ 
&  \leq \int_{\R^2}  \left(  (z_1+z_2)^2 
+ \left(|x-z_2|+|x-z_1|\right)^2
+2|z_1+z_2|\right)\phi(z_1)\phi(z_2)dz_1dz_2 \\ 
&  =2 + \int_{\R^2}  \left(  \left(x-z_2\right)^2+ 
2|x-z_2||x-z_1|+\left(x-z_1\right)^2\right)\phi(z_1)\phi(z_2)dz_1dz_2 
+ \frac{4}{\sqrt{\pi}}\\ &  =2 + 2+2x^2+2\left(\int_{-\infty}^\infty 
|x-z_1|\phi(z_1)dz_1\right)\left(\int_{-\infty}^\infty   
|x-z_2|\phi(z_2)dz_1dz_2\right) + \frac{4}{\sqrt{\pi}}\\ 
&= 4+2x^2+2\left(2x\Phi(x) - x  + 2\phi(x)\right)^2 + 
\frac{4}{\sqrt{\pi}}  \\ &\leq C_1+C_2|x|+C_3x^2. 
\end{align*} 
Hence, 
$$\E\left[h_0(W_t)\right]\leq 
C_1+C_2\E\left[|W_t|\right]+C_3\E\left[W_t^2\right]
=C_1+C_2\sqrt{t}+C_3t. $$
In particular, Condition 5 in Proposition~\ref{prop:ItoTanakaF} holds 
true. 

Similarly, 
$$ |h_1(x)| = \left|2\left(\int_{-\infty}^{\infty}  
|x^2-z_2^2|\phi(z_2)dz_2\right)\phi(x)  \right| 
\leq 2\left(x^2+1\right). $$
Due to the symmetries of the weights 
$h_1=h_2$, and therefore for $i=1,2$ $$\E\left[h_i(W_t)\right]\leq 
C_1+C_2\E\left[W_t^2\right]=C_1+C_2t,$$ and Condition 6 in 
\ref{prop:ItoTanakaF} also holds true. We therefore are in a position 
to use Proposition~\ref{prop:ItoTanakaF} and get
$$\E\bigg[ |G_{i}+G_{i+1}|q(W_{\frac{n-4}{n^2}},G_{i},G_{i+1})
|W_{\frac{n-4}{n^2}}|\bigg] \leq C\frac{n-4}{n^2}, $$ 
from which
$$ E_2 \leq \frac{C(n-4)}{n^3}.$$ 

Similar arguments allow us to show that $ E_i\leq\frac{C(n-4)}{n^3}$ 
for $i=1,3,4$. That ends the proof 
of~\eqref{eq:bound-for-the-rest-q-times-q}.

\bibliographystyle{plain}
\bibliography{Biblio}

\end{document}